\documentclass[pdflatex]{sn-jnl}
\usepackage{graphicx}%
\usepackage{multirow}%
\usepackage{amsmath,amssymb,amsfonts}%
\usepackage{amsthm}%
\usepackage{mathrsfs}%
\usepackage[title]{appendix}%
\usepackage{xcolor}%
\usepackage{textcomp}%
\usepackage{manyfoot}%
\usepackage{booktabs}%
\usepackage{algorithm}%
\usepackage{algorithmicx}%
\usepackage{algpseudocode}%
\usepackage{listings}%
\usepackage[all,pdf]{xy}
\usepackage{scrextend}
\usepackage{fourier}
\deffootnote{0em}{1em}

\theoremstyle{plain}
\newtheorem{theorem}{Theorem}[section]
\newtheorem{lemma}[theorem]{Lemma}
\newtheorem{corollary}[theorem]{Corollary}
\newtheorem{proposition}[theorem]{Proposition}
\theoremstyle{definition}
\newtheorem{definition}[theorem]{Definition}
\newtheorem{example}[theorem]{Example}

\newtheorem{remark}[theorem]{Remark}

\usepackage{tikz-cd}
\usetikzlibrary{cd}

\def\A{\mathcal{A}}

\makeatletter
\g@addto@macro\normalsize{%
  \setlength{\abovedisplayskip}{8pt}
  \setlength{\belowdisplayskip}{8pt}
  \setlength{\abovedisplayshortskip}{4pt}
  \setlength{\belowdisplayshortskip}{4pt}
}
\makeatother

\begin{document}
\normalsize
\title[An Equivariant Coarse Green–Julg Theorem for Proper Actions]{\bfseries An Equivariant Coarse Green–Julg Theorem \break for Proper Actions}

\author[1]{\fnm{Qin} \sur{Wang}}\email{qwang@math.ecnu.edu.cn}

\author*[1]{\fnm{Chen} \sur{Zhang}}\email{52275500018@stu.ecnu.edu.cn}

\affil[1]{%
  \orgdiv{Research Center for Operator Algebras, School of Mathematical Sciences},
  \orgname{East China Normal University},
  \orgaddress{\street{No. 500, Dongchuan Road},\city{Shanghai},\postcode{200241},\country{China}}%
}

\abstract{
We prove an equivariant coarse Green--Julg theorem for proper isometric actions of countable discrete groups on bounded geometry metric spaces. We introduce the $(\Gamma, X)$-controlled proper $C^{*}$-algebra, and show that the evaluation-at-zero map between the equivariant localization and Roe algebras with such coefficients induces an isomorphism on $K$-theory.}

\maketitle
\tableofcontents

\section{Introduction}
The Green--Julg theorem was originally established by P. Julg \cite{Julg71,Jul81} for actions of a compact group $G$ on a $C^{*}$-algebra $\mathcal{A}$. It identifies the equivariant $K$-theory of $\mathcal{A}$ with the $K$-theory of the corresponding crossed product:
\[K_{*}^{G}(\mathcal{A})\cong K_{*}(\mathcal{A}\rtimes{G}).\]
This result was generalized to induced actions through Green's imprimitivity theory \cite{Gre78,Gre80}, and it plays a foundational role in the Baum--Connes conjecture with coefficients \cite{BC88,BCH94}.

Let $\Gamma$ be a countable discrete group and let $\mathcal{A}$ be a $\Gamma$-$C^{*}$-algebra. The Baum--Connes conjecture with coefficients in $\mathcal{A}$ states that the assembly map
\[\mu_{\mathcal{A}}:K_{*}^\Gamma(\underline{E}\Gamma,\mathcal{A})\longrightarrow K_{*}(\mathcal{A}\rtimes_r\Gamma)\]
is an isomorphism.
Here, $\underline{E}\Gamma$ is the classifying space for proper $\Gamma$-actions \cite{BCH94}, and $K_{*}^\Gamma(\underline{E}\Gamma,\mathcal{A})$ is defined as the direct limit of $KK^\Gamma(C_0(X),\mathcal{A})$ taken over all $\Gamma$-compact subsets $X\subseteq\underline{E}\Gamma$.
The notation $KK^\Gamma(C_0(X),\mathcal{A})$ denotes Kasparov's equivariant $KK$-theory \cite{Kas88}, which classifies $\Gamma$-equivariant Kasparov cycles from $C_{0}(X)$ to $\mathcal{A}$. 
In particular, when $\Gamma$ acts properly and cocompactly on a specific model $X$ for $\underline{E}\Gamma$, this $X$ itself is a $\Gamma$-compact subset, and thus the direct limit defining the topological $K$-theory reduces to $KK^\Gamma(C_0(X),\mathcal{A})$. The assembly map takes the form 
$$\mu_{\mathcal{A}}: KK^\Gamma(C_0(X),\mathcal{A})\longrightarrow K_{*}(\mathcal{A}\rtimes_r\Gamma).$$ 
For proper coefficient algebras in the sense of Kasparov, the assembly map is known to be an isomorphism \cite{Kas88,GHT00,HigsonKasparov2001}. This isomorphism is a generalization of the Green--Julg theorem.
Indeed, if we consider a compact group $G$ acting on a single point $pt$, the action is automatically proper and cocompact. Using the fact that $KK^G(C_0(pt),\mathcal{A})$ is isomorphic to $K_{*}^{G}(\mathcal{A})$ \cite{Bla98}, the isomorphism reduces to the classical Green--Julg theorem. If, in addition, $\mathcal{A}=\mathbb{C}$, one further recovers the representation ring isomorphism $R(G)\cong K_{0}(C^{*}_{r}(G))$ \cite{Atiyah69,Seg68}.
The Green--Julg theorem has also been generalized outside the $KK$-theoretic framework. 
In \cite{GHT00}, E. Guentner, N. Higson, and J. Trout established a version of the Green--Julg theorem within the equivariant $E$-theory. Subsequently, building upon the constructions of N. Higson and G. Kasparov, J. Tu \cite{Tu99} generalized the Green--Julg theorem from compact groups to foliations whose holonomy groupoids (and more generally, associated microbundles) are Hausdorff and amenable. We refer the reader to \cite{Par13, Tu04, HG04, CELY17} for further developments in these directions.

In parallel to the group case, G. Yu \cite{Yu95} formulated a coarse analogue of the Baum--Connes conjecture for metric spaces. A key ingredient is the localization algebra $C_{L}^{*}(X)$ introduced by Yu \cite{Yu97}, whose $K$-theory models the locally finite $K$-homology of $X$. To capture the large scale geometry of a metric space $X$, one considers the Rips complexes $P_{d}(X)$ as the scale parameter $d>0$ tends to infinity. The coarse Baum–Connes conjecture asserts that the homomorphism induced by the evaluation-at-zero map
\[e_{*}:\lim_{d\to\infty} K_{*}\bigl(C_{L}^{*}(P_{d}(X))\bigr)\longrightarrow\lim_{d\to\infty} K_{*}\bigl(C^{*}(P_{d}(X))\bigr),\]
is an isomorphism, where $C_{L}^{*}(P_{d}(X))$ and $C^{*}(P_{d}(X))$ denote the localization algebra \cite{Yu97} and Roe algebra \cite{Roe93} associated with $P_{d}(X)$, respectively. A fundamental theorem of G. Yu \cite{Yu00} establishes that the coarse Baum–Connes conjecture holds for discrete metric spaces with bounded geometry that admit a coarse embedding into a Hilbert space.

Combining the perspectives of the coarse Baum--Connes conjecture and the classical Green--Julg theorem leads to the equivariant twisted coarse Baum--Connes conjecture (Cf. \cite{GLWZ23, FZ24, DFW24}). 
Given a countable discrete group $\Gamma$ acting properly on a metric space $X$ and a $\Gamma$-$C^{*}$-algebra $\mathcal{A}$, the conjecture asserts that the homomorphism induced by evaluation-at-zero
\[e_{*}^{\mathcal{A}}:\lim_{d\to\infty} K_{*}\bigl(C_{L}^{*}(P_{d}(X),\mathcal{A})^\Gamma\bigr) \longrightarrow\lim_{d\to\infty} K_{*}\bigl(C^{*}(P_{d}(X),\mathcal{A})^\Gamma\bigr).\]
is an isomorphism.
In the localization algebra approach to coarse Baum--Connes type conjectures \cite{Yu97,Yu00}, the resulting isomorphism is an equivariant coarse analogue of the Green--Julg theorem. 

Recently, L. Guo and Q. Wang \cite{guo2022coarsegeometricapproachgeneralised} gave a coarse geometric proof of the generalized Green--Julg theorem by viewing $\Gamma$ as a metric space with a cocompact left action. Subsequently, L. Guo, Z. Luo, Q. Wang, and Y. Zhang \cite{GLWZ23} introduced the notion of (approximately) coarsely proper algebras and established a coarse geometric analogue of the generalized Green--Julg theorem for metric spaces that admit a coarse embedding into Hilbert space.
In this paper, we generalize these results to proper, not necessarily cocompact, group actions. Our main result is the following theorem.

\begin{theorem}[Main theorem]\label{main theorem}
Let $X$ be a discrete metric space with bounded geometry, and let $\Gamma$ be a countable discrete group acting properly and isometrically on $X$. If the action has equivariant bounded geometry, then for any $(\Gamma, X)$-controlled proper $C^{*}$-algebra $\mathcal{A}$, the equivariant coarse Green--Julg theorem holds, that is, the evaluation map 
\[e_{*}^{\mathcal{A}}:\lim_{d\to\infty}K_{*}\bigl(C^{*}_{L}(P_{d}(X),\mathcal{A})^\Gamma\bigr)\longrightarrow\lim_{d\to\infty} K_{*}\bigl(C^{*}(P_{d}(X),\mathcal{A})^\Gamma\bigr)\]
is an isomorphism.
\end{theorem}

The proof of Theorem \ref{main theorem} relies on several tools. To link the large scale geometry of the base space with the topological structure of the coefficient algebras, we introduce a notion of $(\Gamma, X)$-controlled proper $C^{*}$-algebras. Using the equivariant bounded geometry condition introduced in \cite{GWWY24}, we then reduce the global problem to local models associated with finite stabilizer subgroups. Finally, applying Yu's localization algebra techniques \cite{Yu97}, we prove the isomorphisms for these local models and glue them together via Mayer--Vietoris sequences and the Five Lemma.

\section{Controlled Proper algebras and Equivariant Roe Algebras}

In this section, we set up the geometric and algebraic framework required for our main results.
An essential geometric condition we will use is the notion of equivariant bounded geometry, introduced by L. Guo, Q. Wang, J. Wu, and G. Yu in \cite{GWWY24}. This condition generalizes the classical bounded geometry of a metric space to the equivariant setting and provides the necessary control over the group action in the absence of a compact quotient. 
\begin{definition}[\cite{GWWY24}]\label{def:equiv-bdd-geom}
Let $X$ be a discrete metric space with bounded geometry and let $\Gamma$ be a countable discrete group. A proper isometric action $\Gamma \curvearrowright X$ is said to have \textit{equivariant bounded geometry} if for any $r>0$, there exists an integer $N_r>0$ such that
\[\#\{\gamma\in\Gamma\mid B(\gamma x,r)\cap B(x,r)\neq\emptyset\}<N_{r}\]
for all $x\in X$ (here, $\#A$ denotes the number of elements in a set $A$).
\end{definition}
\begin{remark}\label{rem:equivariant-bounded-geometry}
We make several observations regarding this definition:
\begin{enumerate}
\item[(1)] If $\Gamma$ is the trivial group, the equivariant bounded geometry condition is trivially satisfied. 
\item [(2)] If $X$ is $\Gamma$-compact, then the action $\Gamma \curvearrowright X$ automatically has equivariant bounded geometry. Indeed, choosing a compact subset $K \subseteq X$ with $X=\Gamma\cdot K$, the properness of the action implies that for any $r>0$, only finitely many $\gamma\in\Gamma$ can move a point of $K$ by a distance less than $2r$.
\item [(3)] If the action has equivariant bounded geometry, then there exists $N>0$ such that $\#\Gamma_{x}<N$ for all $x\in X$, where $\Gamma_x=\{\gamma\in\Gamma\mid\gamma x=x\}$ denotes the stabilizer subgroup of $x$. Indeed, fix any $r>0$. If $\gamma\in\Gamma_x$, then $\gamma x=x$, so $B(\gamma x,r)=B(x,r)$ and in particular $B(\gamma x,r)\cap B(x,r)\neq\emptyset$. This gives 
\[\Gamma_x\subseteq\{\gamma \in\Gamma\mid B(\gamma x,r)\cap B(x,r)\neq\emptyset\},\]
and the defining bound $N_r$ immediately yields $\#\Gamma_x<N_r$. Conversely, if the stabilizers have a uniform upper bound $N$ and the space $X$ has bounded geometry, then the action has equivariant bounded geometry. Indeed, for any $R>0$, let $C_{2R}$ be the maximum cardinality of a ball of radius $2R$ in $X$. Since $B(\gamma x,R)\cap B(x,R)\neq\emptyset$ implies $d(\gamma x, x)<2R$, we have
\[\#\{\gamma\in\Gamma\mid B(\gamma x,R)\cap B(x,R)\neq \emptyset\}\leq\#\{\gamma\in\Gamma\mid d(\gamma x, x)<2R\}\leq C_{2R}\cdot N.\]
Here the last inequality holds because the orbit $\Gamma\cdot x$ meets $B(x,2R)$ in at most $C_{2R}$ points, and each lies in the image of at most $\#\Gamma_x\le N$ group elements.
In particular, if $\Gamma$ is torsion free, then since the properness of the action forces all stabilizers to be finite, every stabilizer must be trivial. Thus, the action is free and has equivariant bounded geometry.
\end{enumerate}
\end{remark}
To incorporate coefficients in a geometrically coherent manner, we shall introduce the notion of a $(\Gamma, X)$-controlled proper $C^{*}$-algebra. This notion generalizes the proper $\Gamma$-$C^{*}$-algebra defined in \cite{GLWZ23}. The key idea is to couple the large scale geometry of the base space $X$ with a proper $\Gamma$-space $M$ associated to the coefficient algebra.
\begin{definition}\label{def:proper-Gamma-X-algebra}
Let $\Gamma$ be a countable discrete group acting properly and isometrically on a discrete metric space $X$ with bounded geometry. 
A $\Gamma$-$C^{*}$-algebra $\mathcal{A}$ is called a \emph{$(\Gamma, X)$-controlled proper $C^{*}$-algebra} if there exists a second countable, locally compact, Hausdorff space $M$ with a proper $\Gamma$-action, and non-decreasing functions
$\rho_{\pm}\colon[0,\infty)\to[0,\infty)$ with $\rho_-(r)>0$ for all $r>0$
and $\lim_{r\to\infty}\rho_-(r)=\infty$, such that:
\begin{enumerate}
\item[(1)] there is a $\Gamma$-equivariant $*$-homomorphism $\varphi$ from $C_0(M)$ to $Z(\mathcal{M}(\mathcal{A}))$,
the center of the multiplier 
algebra $\mathcal{M}(\mathcal{A})$, such that
$\overline{\varphi(C_0(M))\cdot\mathcal{A}}=\mathcal{A}$;
\item[(2)]there is a family of open subsets $\{O_{x,r}\}_{x\in X,r>0}$ of $M$ with $O_{x,r}\subseteq O_{x,s}$ for $r\le s$, such that:
\begin{enumerate}
\item[(a)] for every $x\in X$, $\bigcup_{r>0}O_{x,r}=M$;
\item[(b)] if $r>0$ and $d(x_1,x_2)>2r$, then $O_{x_1,\rho_-(r)}\cap O_{x_2,\rho_-(r)}=\emptyset$;
\item[(c)] for all $r>0$ and $x_1,x_2\in X$, $O_{x_1,r}\subseteq O_{x_2,\rho_+(d(x_1,x_2))+r}$;
\item[(d)] $\gamma O_{x,r}=O_{\gamma x,r}$ for all $\gamma\in\Gamma$, $x\in X$, $r>0$.
\end{enumerate}
\end{enumerate}
\end{definition}

\begin{example}\label{ex:Hilbert-Space-Algebra}
We present an example of a $(\Gamma,X)$-controlled proper $C^*$-algebra, following the construction of the $C^*$-algebra $\mathcal{A}(M)$ in \cite{GWWY24} (see also \cite{GHT00,GWY21}), specialized to the case where $M$ is a Hilbert space.\par
Let $H_0$ be a separable infinite-dimensional real Hilbert space on which $\Gamma$ acts by isometries, and let $f_0\colon X\to H_0$ be a $\Gamma$-equivariant coarse embedding. That is, $f_{0}(\gamma x)=\gamma\cdot f_0(x)$ for all $\gamma\in\Gamma$, and there exist non-decreasing control functions $\rho_\pm^{0}\colon [0,\infty)\to[0,\infty)$ with $\rho^{0}_{-}(s)\to\infty$ as $s\to\infty$, such that
\[
\rho_{-}^{0}(d(x,y))\le\|f_0(x)-f_0(y)\|\le\rho_{+}^{0}(d(x,y))
\]
for $x,y\in X$.
By the Mazur--Ulam theorem \cite{mazur1932sur}, this action is affine. Thus, there exist an orthogonal representation $\pi_0\colon\Gamma\to O(H_0)$ and a $1$-cocycle $b\colon\Gamma\to H_{0}$ such that 
$\gamma\cdot v=\pi_{0}(\gamma)v+b(\gamma)$
for all $\gamma\in\Gamma$ and $v\in H_{0}$.
Since $X$ has bounded geometry, $X$ is countable. Define $$H:=H_{0}\oplus\ell^2(X;\mathbb R).$$
Let $\lambda_{X}\colon\Gamma\to O(\ell^2(X;\mathbb R))$ be the representation induced by the action of $\Gamma$ on $X$.
Define $\pi:=\pi_0\oplus\lambda_X$ and $c(\gamma):=(b(\gamma),0)$. We equip $H$ with the affine isometric action
\[
\gamma\cdot v:=\pi(\gamma)v+c(\gamma),
\ \text{where}\ \gamma\in\Gamma,\ v\in H.
\]
Define $f(x)=(f_0(x),\delta_x)$. 
Then $f$ is $\Gamma$-equivariant, the action on $H$ is affine isometric, 
and $$\|f(x)-f(y)\|^2=\|f_0(x)-f_0(y)\|^2+2\ge 2$$ whenever $x\neq y$. Define
\[
\rho_{+}(s)=\sup\{\|f(x)-f(y)\| \mid d(x,y)\le s\},\quad
\rho_{f,-}(s)=\min\bigl\{s,\ \inf\{\|f(x)-f(y)\| \mid d(x,y)\ge s\}\bigr\}.
\]
Both are non-decreasing, $\lim_{s\to\infty}\rho_{f,-}(s)=\infty$, and $\rho_{f,-}(d(x,y))\le\|f(x)-f(y)\|\le\rho_+(d(x,y))$.

Fix $x_0\in X$. Since $\Gamma\curvearrowright X$ is proper, $\{\gamma\in\Gamma\mid d(\gamma x_0,x_0)\le s\}$ is finite for every $s$. Note that
\[
\|\gamma f(x_0)-f(x_0)\|=\|f(\gamma x_0)-f(x_0)\|\ge\rho_{f,-}(d(\gamma x_0,x_0)).
\]
Then $\{\gamma\in\Gamma\mid\|\gamma f(x_0)-f(x_0)\|\le R\}$ is finite for every $R$. As $\gamma$ is an isometry, $\|\gamma v-v\|\ge\|\gamma f(x_0)-f(x_0)\|-2\|v-f(x_0)\|$ for all $v\in H$. Hence, $\{\gamma\in\Gamma\mid\gamma B(v,R)\cap B(v,R)\ne\emptyset\}$ is finite for all $v,R$, which implies $\Gamma\curvearrowright H$ is a proper action.

For $(v,t)\in H\times[0,\infty)$, set $H_{t}:=H\oplus\mathbb{R}_{t}$, where $\mathbb{R}_{t}=\mathbb{R}$ if $t>0$ and $\{0\}$ if $t=0$. Define
\[
\Pi_b(H)=\Bigl\{\sigma\in\prod_{(v,t)\in H\times[0,\infty)}\operatorname{Cliff}_{\mathbb C}(H_t)\ \mid\ \sup\|\sigma(v,t)\|<\infty\Bigr\}.
\]
For $w\in H$, define the Clifford generator: $$C_{w}(v,t)=(w-v,\,t)\in H_{t}\subset\operatorname{Cliff}_{\mathbb C}(H_t).$$
Let $\mathcal{S}=C_0(\mathbb{R})$, equipped with the standard $\mathbb{Z}_{2}$-grading. By the functional calculus, the Bott map
\[
\beta_w\colon\mathcal S\to\Pi_b(H),\qquad \beta_w(f)(v,t)=f\bigl(C_w(v,t)\bigr)
\]
is a graded $*$-homomorphism. 
Now, we define $\mathcal{A}(H)$ to be the $C^{*}$-subalgebra of $\Pi_b(H)$ generated by $\{\beta_w(f)\mid w\in H,\ f\in\mathcal S\}$. 
For $\gamma\in\Gamma$ and $\sigma\in\Pi_b(H)$, define
\[
(\gamma_{*}\sigma)(v,t)=\operatorname{Cliff}_{\mathbb{C}}\bigl(\pi(\gamma)\oplus 1_{\mathbb R_{t}}\bigr)\bigl(\sigma(\gamma^{-1}v,t)\bigr).
\]
This is an action by automorphisms. Since $\pi(\gamma)\bigl(w-\gamma^{-1}v\bigr)=\gamma w-v$, we get $\gamma_{*}\beta_w(f)=\beta_{\gamma w}(f)$, so $\gamma_{*}$ preserves $\mathcal{A}(H)$. As $\Gamma$ is discrete, $\mathcal{A}(H)$ is a $\Gamma$-$C^{*}$-algebra.

Equip $Y=H\times[0,\infty)$ with the weakest topology in which all functions
\[
\varphi_{w,g}(v,t)=g\bigl(\sqrt{\|v-w\|^2+t^2}\bigr),\qquad w\in H,\ g\in C_0([0,\infty)),
\]
are continuous. Since
\[
\|v-w\|^2+t^2-(\|v\|^2+t^2)=\|w\|^2-2\langle v,w\rangle,
\]
this is the weakest topology for which $(v,t)\mapsto\langle v,w\rangle$ ($w\in H$) and $(v,t)\mapsto \|v\|^2+t^2$ are continuous. Thus, $Y$ is a second countable, locally compact Hausdorff space.
The relation $\varphi_{w,g}\circ\gamma=\varphi_{\gamma^{-1}\cdot w,g}$ guarantees that $\Gamma$ acts on $Y$ by homeomorphisms via $\gamma\cdot(v,t)=(\gamma\cdot v,t)$. For any compact subset $K\subset Y$, $K$ is contained in $\{(v,t)\in Y\mid\|v\|^2+t^2\le R^2\}$ for some $R>0$. If $\gamma K\cap K\neq\emptyset$, there exists $(v,t)\in K$ with $\|v\|\le R$ and $\|\gamma\cdot v\|\le R$. We have that
\[
\|\gamma\cdot 0-0\|\le\|\gamma\cdot 0-\gamma\cdot v\|+\|\gamma\cdot v\|\le 2R.
\]
Since the action on $H$ is proper, the properness of the action on $Y$ follows immediately.

Observe that each $\varphi_{w,g}$ vanishes at infinity. Since this family separates the points of $Y$ and is nowhere vanishing, the Stone--Weierstrass theorem ensures that the $C^{*}$-algebra it generates is exactly $C_{0}(Y)$.
Note that
\[
\varphi_{w,g}=\beta_w(\tilde g),\ \text{with}\ \tilde g(u)=g(|u|)\in\mathcal{S}_{\mathrm{ev}},
\]
a scalar-valued element of $\mathcal{A}(H)$. Hence $C_0(Y)\subset\mathcal{A}(H)\subset\Pi_b(H)$. For $\varphi\in C_{0}(Y)$ and $a\in\mathcal{A}(H)$ we have $\varphi a\in\mathcal{A}(H)$ and $\varphi a=a\varphi$, since $\varphi$ takes values in $\mathbb C\cdot1$. So multiplication by $\varphi$ defines a central multiplier, giving a $*$-homomorphism $C_{0}(Y)\to Z(\mathcal{M}(\mathcal{A}(H)))$. It is $\Gamma$-equivariant since
\[(\gamma_{*}\varphi)(v,t)=\varphi(\gamma^{-1}\cdot v,t),
\]
as $\operatorname{Cliff}\bigl(\pi(\gamma)\oplus 1\bigr)$ is the identity on scalars.
Take $e_n\in C_c(\mathbb R)_{\mathrm{ev}}$ with $0\le e_n\le1$ and $e_n\to 1$ uniformly on compacts. For $f\in\mathcal S$,
\[\beta_w(f)\beta_w(e_n)=\beta_w(fe_n)\to\beta_w(f),
\]
and $\beta_w(e_n)\in C_0(Y)$. Hence every generator lies in $\overline{C_0(Y)\cdot\mathcal{A}(H)}$, which is therefore all of $\mathcal{A}(H)$.

Put $\rho_-(r)=\tfrac12\rho_{f,-}(2r)$ and $\rho_+$ as above. For $x\in X$, $r>0$, set
\[O_{x,r}=\bigl\{(v,t)\in Y\mid\|v-f(x)\|^2+t^{2}<r^{2}\bigr\}.
\]
Condition (2)(a) is clear.
For (2)(b), let $d(x_1,x_2)>2r$ and suppose $(v,t)\in O_{x_1,\rho_-(r)}\cap O_{x_2,\rho_-(r)}$. Then
\[
\|f(x_1)-f(x_2)\|\le\|f(x_1)-v\|+\|v-f(x_2)\|<2\rho_-(r)=\rho_{f,-}(2r),
\]
whereas $\|f(x_1)-f(x_2)\|\ge\rho_{f,-}(d(x_1,x_2))\ge\rho_{f,-}(2r)$, a contradiction.
For (2)(c), for $(v,t)\in O_{x_1,r}$,
\[
\sqrt{\|v-f(x_2)\|^2+t^2}\le\sqrt{\|v-f(x_1)\|^2+t^2}+\|f(x_1)-f(x_2)\|<r+\rho_{+}(d(x_1,x_2)).
\]
For (2)(d), since $\gamma$ is an isometry of $H$ and $f(\gamma x)=\gamma\cdot f(x)$,
\[
\gamma O_{x,r}=\{(w,t)\mid\|\gamma^{-1}\cdot w-f(x)\|^2+t^2<r^2\}=\{(w,t)\mid \|w-f(\gamma x)\|^2+t^2<r^2\}=O_{\gamma x,r}.
\]
Hence, $\mathcal{A}(H)$ is a $(\Gamma,X)$-controlled proper $C^*$-algebra as defined in Definition \ref{def:proper-Gamma-X-algebra}.
\end{example}

Now, we are ready to construct the equivariant Roe algebra with coefficients in a
$(\Gamma,X)$-controlled proper $C^*$-algebra $\mathcal{A}$.
Let $H_0$ be a fixed infinite-dimensional separable Hilbert space,
equipped with the trivial $\Gamma$-action. Define
\[H_{X,\mathcal{A}}:=
\ell^2(X)\otimes\ell^2(\Gamma)\otimes H_0\otimes\mathcal{A},\]which is the completion of the algebraic tensor product with respect to the $\mathcal{A}$-valued inner product:
\[\left\langle
\delta_x\otimes\delta_\gamma\otimes v\otimes a,\,
\delta_y\otimes\delta_{\gamma'}\otimes w\otimes b
\right\rangle_{\mathcal{A}}=
\delta_{x,y}\delta_{\gamma,\gamma'}
\langle v,w\rangle_{H_0}\,a^*b,\]
and is equipped with the right action:
\[(\delta_x\otimes\delta_\gamma\otimes v\otimes a)b=
\delta_x\otimes\delta_\gamma\otimes v\otimes ab.\]
This makes $H_{X,\mathcal{A}}$ a right Hilbert $\mathcal{A}$-module.
Here $x,y\in X$, $\gamma,\gamma'\in\Gamma$, $v,w\in H_0$,
$a,b\in\mathcal{A}$, and $\delta_{x,y}$ and
$\delta_{\gamma,\gamma'}$ denote the Kronecker symbols.
The diagonal $\Gamma$-action on $H_{X,\mathcal{A}}$ is defined by
\[U_\eta(\delta_x\otimes\delta_\gamma\otimes v\otimes a)
=\delta_{\eta x}\otimes\delta_{\eta\gamma}
\otimes v\otimes(\eta\cdot a),
\ \eta\in\Gamma.\]
It is compatible with the Hilbert $\mathcal{A}$-module structure:
\[U_\eta(\xi a)=(U_\eta\xi)(\eta\cdot a),
\qquad\langle U_\eta\xi,U_\eta\zeta\rangle_{\mathcal{A}}
=\eta\cdot\langle\xi,\zeta\rangle_{\mathcal{A}}.\]
Thus, $H_{X,\mathcal{A}}$ is a $\Gamma$-equivariant Hilbert
$\mathcal{A}$-module.

For Hilbert $\mathcal{A}$-modules $E$ and $F$, let
$\mathcal{L}(E,F)$ denote the space of adjointable
$\mathcal{A}$-linear operators from $E$ to $F$, and 
$\mathcal{L}(E):=\mathcal{L}(E,E)$.
There is a faithful, non-degenerate, equivariant $*$-representation $\psi$ from $C_0(X)$ to
$\mathcal{L}(H_{X,\mathcal{A}})$
defined by pointwise multiplication:
\[\psi(f)(\delta_x\otimes\delta_\gamma\otimes v\otimes a)
=f(x)\,\delta_x\otimes\delta_\gamma\otimes v\otimes a.\]
The equivariance is given by
$U_\eta\psi(f)U_\eta^{-1}=\psi(\eta\cdot f)$ and 
$(\eta\cdot f)(x)=f(\eta^{-1}x)$.
Denote by $\mathcal K(H_{X,\mathcal{A}})$ the algebra of compact
adjointable operators on $H_{X,\mathcal{A}}$. 
We write
\[H_{X,\mathcal{A}}=
\bigoplus_{x\in X}H_x,\ \text{where}\
H_x:=\mathbb C\delta_x\otimes\ell^2(\Gamma)
\otimes H_0\otimes\mathcal{A}.\]
Let
$P_x:=\psi(\chi_{\{x\}})$
be the projection onto $H_x$. For
$T\in\mathcal{L}(H_{X,\mathcal{A}})$, define its matrix coefficient by
$$T(x,y):=\left.P_xTP_y\right|_{H_y}
\in\mathcal{L}(H_y,H_x),\ x,y\in X.$$ 
Then $T$ can be represented by the matrix
$\bigl(T(x,y)\bigr)_{x,y\in X}$.
Since $C_0(M)$ acts on $\mathcal{A}$ by central multipliers,
each $f\in C_0(M)$ defines an adjointable operator $\phi(f)$ on
$H_{X,\mathcal{A}}$ by
\[\phi(f)(\delta_x\otimes\delta_\gamma\otimes v\otimes a)
:=\delta_x\otimes\delta_\gamma\otimes v\otimes fa.\]
For $x,y\in X$, the support of $T(x,y)$ is the closed subset
\[\operatorname{Supp}\bigl(T(x,y)\bigr):=M\setminus\bigl\{m\in M\mid
\phi(f)\,T(x,y)=0\ \text{for some}\ f\in C_0(M)\ \text{with}\ f(m)\neq0\bigr\}.\]

\begin{definition}\label{def:equiv-alg-roe-X}
The equivariant algebraic Roe algebra with coefficients in a $(\Gamma,X)$-controlled proper $C^{*}$-algebra $\mathcal{A}$, denoted by $\mathbb{C}[X,\mathcal{A}]^{\Gamma}$, is defined as the $*$-algebra of all bounded functions $T:X\times X\to \mathcal{K}(\ell^2(\Gamma)\otimes H_0)\otimes \mathcal{A}$ such that:
\begin{enumerate}
\item[(1)]for any bounded subset $B\subset X$, the set
$\{(x,y)\in(B\times B)\mid T(x,y)\neq 0 \}$
is finite;
\item[(2)]there exists $L>0$ such that 
$\#\{y\in X\mid T(x,y)\neq 0\}<L$ and
$\#\{y\in X\mid T(y,x)\neq 0\}<L$
for all $x\in X$;
\item[(3)] there exists $R>0$ such that $T(x,y)=0$ whenever $d(x,y)>R$ for any $x,y\in X$ (the minimal $R$ is called the \emph{propagation} of $T$, denoted by $\operatorname{prop}(T)$).
\item[(4)]$T$ is $\Gamma$-equivariant, i.e., $T(\gamma x,\gamma y)=\gamma\cdot T(x, y)$ for all $\gamma\in \Gamma$ and $x, y\in X$, where $\gamma\cdot$ denotes the diagonal $\Gamma$-action on $\mathcal{K}(\ell^2(\Gamma)\otimes H_0) \otimes\mathcal{A}$.
\item[(5)]there exists $r>0$ such that $\operatorname{Supp}(T(x, y))\subseteq O_{x,r}$ for all $x, y\in X$, where $O_{x,r}$ is given in Definition \ref{def:proper-Gamma-X-algebra}.
\end{enumerate}
\end{definition}

\begin{remark}\label{rem:closure-algebra}
Condition (5) ensures that the coefficient part reflects the coarse geometry of $X$. We briefly verify that
$\mathbb{C}[X,\mathcal{A}]^\Gamma$ is a $*$-algebra.
Let $S,T\in\mathbb{C}[X,\mathcal{A}]^\Gamma$, with propagation bounds
$R_S,R_T$ and coefficient support bounds $r_S,r_T$. For the adjoint, we have
$T^{*}(x,y)=T(y,x)^{*}$,
and hence, whenever $d(x,y)\leq R_T$,
$\operatorname{Supp}(T^*(x,y))
\subseteq O_{y,r_T}
\subseteq O_{x,\rho_+(R_T)+r_T}$
(by condition (2)(c) of
Definition \ref{def:proper-Gamma-X-algebra}). Thus, $T^{*}$ satisfies
(5).
For the product,
$(ST)(x,y)=\sum_{z\in X}S(x,z)T(z,y)$.
The sum is finite, and
$\operatorname{prop}(ST)\leq R_S+R_T$.
Moreover,
$\operatorname{Supp}\bigl(S(x,z)T(z,y)\bigr)
\subseteq
\operatorname{Supp}(S(x,z))
\subseteq O_{x,r_S}$,
so $ST$ satisfies  (5). The uniform row and column
bounds, local finiteness, and equivariance are preserved under
addition, multiplication, and involution. So,
$\mathbb{C}[X,\mathcal{A}]^\Gamma$ is a well-defined $*$-algebra.
\end{remark}
\begin{definition}\label{def:equiv-roe-X}
The \emph{equivariant Roe algebra with coefficients in a $(\Gamma, X)$-controlled proper $C^{*}$-algebra $\mathcal{A}$}, denoted by $C^{*}(X,\mathcal{A})^{\Gamma}$, is defined as the operator norm completion of $\mathbb{C}[X,\mathcal{A}]^{\Gamma}$ in $\mathcal{L}(H_{X,\mathcal{A}})$, where $\mathcal{L}(H_{X,\mathcal{A}})$ denotes the $C^{*}$-algebra of all adjointable operators on the Hilbert $\mathcal{A}$-module $H_{X,\mathcal{A}}$.
\end{definition}
To formulate the evaluation map in Theorem \ref{main theorem}, we need the corresponding localization algebra. 
\begin{definition}\label{def:equiv-alg-loc-X}
The \emph{algebraic localization algebra with coefficients in a $(\Gamma, X)$-controlled proper $C^{*}$-algebra $\mathcal{A}$}, denoted by $\mathbb{C}_L[X,\mathcal{A}]^{\Gamma}$, is defined as the $*$-algebra of all bounded and uniformly continuous functions $g\colon [0,\infty)\to\mathbb{C}[X,\mathcal{A}]^{\Gamma}$ such that:
\begin{enumerate}
\item[(1)] $\text{sup}_{t\geq 0}\text{prop}(g(t))<\infty$ and $\operatorname{prop}(g(t))\to 0$ as $t\to \infty$;
\item[(2)] there exists $r>0$ (independent of $t$) such that for all $t\in [0,\infty)$ and all $x, y\in X$, $$\operatorname{Supp}(g(t)(x,y))\subseteq O_{x,r},$$ where $O_{x,r}\subseteq M$ is the open set associated to $x\in X$ and $r>0$ as in Definition \ref{def:proper-Gamma-X-algebra};
\end{enumerate}
The \emph{equivariant localization algebra with coefficients} in a $(\Gamma, X)$-controlled proper $C^{*}$-algebra $\mathcal{A}$, denoted by $C^{*}_L(X,\mathcal{A})^{\Gamma}$, is defined as the completion of $\mathbb{C}_L[X,\mathcal{A}]^{\Gamma}$ with respect to the norm
\[\|g\|:=\sup_{t\in [0,\infty)}\|g(t)\|_{C^{*}(X,\mathcal{A})^{\Gamma}}.\]
\end{definition}
Let $P_{d}(X)$ be the Rips complex of $X$ at scale $d>0$,
equipped with the semi-spherical metric $d_{P_{d}}$ (Cf.\cite{WY2020}). 
For every $d>0$,
choose a $\Gamma$-invariant countable dense subset
$Z_{d}\subseteq P_{d}(X)$
such that
$Z_{d}\subseteq Z_{d'}$ whenever $d\leq d'$.
For any $x\in P_{d}(X)$, let $V_x\subset X$ denote the finite set of vertices of the unique open simplex containing $x$. 
The map $z\mapsto V_z$ is $\Gamma$-equivariant:
$V_{\gamma z}=\gamma V_z$ for all
$\gamma\in\Gamma$. For $z\in Z_{d}$ and $r>0$, we define
\[O_{z,r}:=\bigcup_{v\in V_z}O_{v,r}\subseteq M.
\]
The sets $O_{z,r}$ are open, increase with $r$, and exhaust $M$.
They are also $\Gamma$-equivariant:
\[O_{\gamma z,r}
=\bigcup_{v'\in V_{\gamma z}}O_{v',r}
 =\bigcup_{v\in V_z}O_{\gamma v,r} 
=\gamma\cdot\bigcup_{v\in V_z}O_{v,r}
 =\gamma O_{z,r}.\]
Furthermore, for any $v_z\in V_z$, condition
\textup{(2)(c)} of
Definition \ref{def:proper-Gamma-X-algebra} gives
$$O_{v_z,r}
\subseteq O_{z,r}
\subseteq O_{v_z,\rho_+(d)+r}.$$
Thus, for each fixed $d$, the collection
$\{O_{z,r}\}_{z\in Z_{d},r>0}$ is uniformly controlled by the
original family $\{O_{x,r}\}_{x\in X,r>0}$.
Using $Z_{d}$ as an indexing set, define the equivariant Hilbert
$\mathcal{A}$-module
\[H_{P_{d}(X),\mathcal{A}}
:=
\ell^{2}(Z_{d})\otimes\ell^2(\Gamma)\otimes H_{0}\otimes\mathcal{A}.
\]
We define the equivariant algebraic Roe algebra
$\mathbb{C}[P_{d}(X),\mathcal{A}]^\Gamma$
via $Z_{d}\times Z_{d}$-matrices satisfying Definition \ref{def:equiv-alg-roe-X}, using the metric $d_{P_{d}}$ and requring $\operatorname{Supp}(T(z,w))\subseteq O_{z,r}$.
Its completion in
$\mathcal{L}(H_{P_{d}(X),\mathcal{A}})$ is denoted by
$C^{*}(P_{d}(X),\mathcal{A})^\Gamma$.
The localization algebra
$C_{L}^{*}(P_{d}(X),\mathcal{A})^\Gamma$
is defined analogously. Then there is a evaluation-at-zero map $$e^{\A}:C^{*}_L(P_{d}(X),\mathcal{A})^\Gamma \longrightarrow C^{*}(P_{d}(X),\mathcal{A})^\Gamma,\quad e^{\A}(g)=g(0).$$
Passing to $K$-theory and taking the direct limit, it induces
\begin{equation}\label{(1)(1)}
e_{*}^{\A}:\lim_{d\to\infty}K_{*}\bigl(C^{*}_L(P_{d}(X),\mathcal{A})^{\Gamma}\bigr) \longrightarrow\lim_{d\to\infty}K_{*}\bigl(C^{*}(P_{d}(X),\mathcal{A})^{\Gamma}\bigr).
\end{equation}

\section{Ideal Decompositions and Slice Separated Families}

To establish the equivariant coarse Green--Julg theorem for proper group actions, now it suffices to prove that map (\ref{(1)(1)}) is an isomorphism.
Following \cite{SW07}, we recall the notion of support. For $T\in\mathbb{C}[P_{d}(X),\mathcal{A}]^{\Gamma}$, its support is the subset of $Z_{d}\times Z_{d}\times M$ given by
\[
\operatorname{Supp}(T):=\{(x,y,m)\in Z_{d}\times Z_{d}\times M\mid m \in\operatorname{Supp}(T(x, y))\}.
\]
For $g\in\mathbb{C}_L[P_{d}(X),\mathcal{A}]^{\Gamma}$, its support is given by
\[
\operatorname{Supp}(g):=\bigcup_{t\in [0,\infty)}\operatorname{Supp}(g(t)).
\]
Let $O$ be a $\Gamma$-invariant open subset of $M$. Define
\[
\mathbb{C}[P_{d}(X),\mathcal{A}]_{O}^{\Gamma}:=\big\{T \in\mathbb{C}[P_{d}(X),\mathcal{A}]^{\Gamma}\mid \operatorname{Supp}(T(x,y))\subseteq O \ \text{for all}x,y\in Z_{d}\big\}.\]
Let $C^{*}(P_{d}(X),\mathcal{A})_O^{\Gamma}$ be the norm closure of $\mathbb{C}[P_{d}(X),\mathcal{A}]_O^{\Gamma}$ in $C^{*}(P_{d}(X),\mathcal{A})^{\Gamma}$.
Similarly, define
\[
\mathbb{C}_L[P_{d}(X),\mathcal{A}]_O^{\Gamma}:=\big\{g \in\mathbb{C}_L[P_{d}(X),\mathcal{A}]^{\Gamma}\mid\operatorname{Supp}(g)\subseteq Z_{d}\times Z_{d}\times O\big\},
\]
and let $C^{*}_L(P_{d}(X),\mathcal{A})_O^{\Gamma}$ be closure of $\mathbb{C}_L[P_{d}(X),\mathcal{A}]_O^{\Gamma}$ in $C^{*}_L(P_{d}(X),\mathcal{A})^{\Gamma}$ with respect to the norm
\[
\|g\|:=\sup_{t\in [0,\infty)}\|g(t)\|_{C^{*}(P_{d}(X),\mathcal{A})^{\Gamma}}.
\]
Since $C_0(M)$ acts by central multipliers, the support condition is preserved under taking adjoints and under multiplication. Consequently, $\mathbb{C}[P_{d}(X),\mathcal{A}]_O^{\Gamma}$ is a two sided $*$-ideal in $\mathbb{C}[P_{d}(X),\mathcal{A}]^{\Gamma}$, and by the norm continuity of the operations, its closure $C^{*}(P_{d}(X),\mathcal{A})_O^{\Gamma}$ is a closed two sided $*$-ideal (similar for localization algebras).
The evaluation-at-zero map then restricts to the $*$-homomorphism:
\[
e:C^{*}_L(P_{d}(X),\mathcal{A})_O^{\Gamma}\longrightarrow C^{*}(P_{d}(X),\mathcal{A})_O^{\Gamma},\quad e(g)=g(0).
\]

\begin{lemma}\label{lem:mv-ideals}
For any two $\Gamma$-invariant open subsets $O_1, O_2 \subseteq M$, we have that
\begin{enumerate}
\item[(1)] $C^{*}(P_{d}(X),\mathcal{A})^{\Gamma}_{O_{1}}+C^{*}(P_{d}(X),\mathcal{A})^{\Gamma}_{O_{2}}=C^{*}(P_{d}(X),\mathcal{A})^{\Gamma}_{O_{1}\cup O_{2}}$;
\item[(2)] $C^{*}(P_{d}(X),\mathcal{A})^{\Gamma}_{O_{1}}\cap C^{*}(P_{d}(X),\mathcal{A})^{\Gamma}_{O_{2}}=C^{*}(P_{d}(X),\mathcal{A})^{\Gamma}_{O_{1}\cap O_{2}}$;
\item[(3)] $C^{*}_{L}(P_{d}(X),\mathcal{A})^{\Gamma}_{O_{1}}+C^{*}_{L}(P_{d}(X),\mathcal{A})^{\Gamma}_{O_{2}}=C^{*}_{L}(P_{d}(X),\mathcal{A})^{\Gamma}_{O_{1}\cup O_{2}}$;
\item[(4)] $C^{*}_{L}(P_{d}(X),\mathcal{A})^{\Gamma}_{O_{1}}\cap C^{*}_{L}(P_{d}(X),\mathcal{A})^{\Gamma}_{O_{2}}=C^{*}_{L}(P_{d}(X),\mathcal{A})^{\Gamma}_{O_{1}\cap O_{2}}$.
\end{enumerate}
\end{lemma}
\begin{proof}[\textbf{Proof.}]\normalsize
For (1), the inclusion $\subseteq$ is clear. For the reverse on the dense algebraic subalgebras,
\[\mathbb{C}[P_{d}(X),\mathcal{A}]^{\Gamma}_{O_{1}\cup O_{2}}\subseteq\mathbb{C}[P_{d}(X),\mathcal{A}]^{\Gamma}_{O_{1}} + \mathbb{C}[P_{d}(X),\mathcal{A}]^{\Gamma}_{O_{2}},\]
we let $T\in\mathbb{C}[P_{d}(X),\mathcal{A}]_{O_1 \cup O_2}^\Gamma$. Since the action $\Gamma\curvearrowright M$ is proper and $O_1,O_2$ are $\Gamma$-invariant open subsets, there exists a $\Gamma$-invariant partition of unity $\{\varphi_1,\varphi_2\}$ on $O_1 \cup O_2$ subordinate to the cover $\{O_1, O_2\}$, with $0\le\varphi_i\le 1$, $\operatorname{Supp}(\varphi_i) \subseteq O_i$ for $i=1, 2$, and $\varphi_1+\varphi_2\equiv 1$ on $O_1\cup O_2$.

Define $T_1,T_2: Z_{d}\times Z_{d}\to\mathcal{K}(\ell^2(\Gamma)\otimes H_0)\otimes\mathcal{A}$ by
\[
T_1(x, y)=(1\otimes\varphi_1)\cdot T(x, y),\quad 
T_2(x, y)=(1\otimes\varphi_2)\cdot T(x, y),\]
where $\varphi_i$ acts as a central multiplier on $\mathcal{A}$ through $C_0(O_1\cup O_2)\subseteq C_0(M)\hookrightarrow Z(M(\mathcal{A}))$, and trivially on the $\mathcal{K}(\ell^2(\Gamma)\otimes H_0)$. Note that:
\begin{itemize}
\item[(a)] $\operatorname{Supp}(T_1(x,y))\subseteq\operatorname{Supp}(\varphi_1)\cap\operatorname{Supp}(T(x,y))\subseteq O_1$;
\item[(b)] $\operatorname{prop}(T_1)\leq \operatorname{prop}(T)<\infty$ since $T_1$ is given by pointwise multiplication;
\item[(c)] For any $\gamma\in\Gamma$, the $\Gamma$-equivariance of $T$ and $\Gamma$-invariance of $\varphi_1$ yields
\[T_1(\gamma x,\gamma y)=(1\otimes\varphi_1)\cdot T(\gamma x,\gamma y)=(1\otimes\varphi_1)\cdot(\gamma \cdot T(x,y))=\gamma\cdot((1\otimes\varphi_1)\cdot T(x, y))=\gamma\cdot T_1(x, y);\]
\item[(d)] Local finiteness and uniform boundedness are inherited from $T$.
\end{itemize}
So $T_1\in\mathbb{C}[P_{d}(X),\mathcal{A}]_{O_1}^\Gamma$. Similarly, $T_2\in \mathbb{C}[P_{d}(X),\mathcal{A}]_{O_2}^\Gamma$. Since $\varphi_1+\varphi_2\equiv 1$ on $O_1\cup O_2$, we have $T=T_1+T_2$. Passing to the norm completions, $C^{*}(P_{d}(X),\mathcal{A})_{O_1}^\Gamma$ and $C^{*}(P_{d}(X),\mathcal{A})_{O_2}^\Gamma$ are closed two sided ideals in $C^{*}(P_{d}(X),\mathcal{A})_{O_1\cup O_2}^\Gamma$. Since the sum of two closed ideals is still closed, the $C^*$-algebra decomposes accordingly.

For (2), the inclusion $\supseteq$ is obvious. For any $T_1 \in\mathbb{C}[P_{d}(X),\mathcal{A}]_{O_1}^\Gamma$ and $T_2\in\mathbb{C}[P_{d}(X),\mathcal{A}]_{O_2}^\Gamma$, their product $T_{1}T_2$ satisfies
\[
\operatorname{Supp}((T_1 T_2)(x, y))\subseteq\bigcup_{z\in Z_{d}}\bigl(\operatorname{Supp}(T_1(x, z))\cap\operatorname{Supp}(T_2(z, y))\bigr)\subseteq O_1\cap O_2.
\]
Thus, $T_1 T_2\in\mathbb{C}[P_{d}(X),\mathcal{A}]_{O_1\cap O_2}^\Gamma$. From the fact that in a $C^*$-algebra, the intersection of two closed two sided ideals coincides with the closed linear span of their product, we obtain $$C^{*}(P_{d}(X),\mathcal{A})_{O_1}^\Gamma\cap C^{*}(P_{d}(X),\mathcal{A})_{O_2}^\Gamma\subseteq C^{*}(P_{d}(X),\mathcal{A})_{O_1\cap O_2}^\Gamma.$$

Case (3) follows by applying the above decomposition pointwise in $t\in[0,\infty)$, and case (4) by the same argument in $C^{*}_{L}(P_{d}(X),\mathcal{A})^\Gamma$.
\end{proof}

Lemma \ref{lem:mv-ideals} reduces Mayer--Vietoris arguments to decompositions of $M$ into $\Gamma$-invariant open subsets. We obtain such subsets from $\Gamma$-equivariant families of open sets indexed by points of a $\Gamma$-invariant subset $X_0\subseteq X$, by taking unions along $\Gamma$-orbits. The families we use are the following.

\begin{definition}\label{definition 3.10}
Let $r>0$ and $X_0$ be a $\Gamma$-invariant countable subset of $X$. A family of open subsets $\{U_x\}_{x\in X_0}$ of $M$ is said to be \emph{$(\Gamma, X_0,r)$-separated} if:
\begin{enumerate}
    \item[(1)] The family $\{U_x\}_{x\in X_0}$ is $\Gamma$-equivariant, i.e., $\gamma U_x=U_{\gamma x}$ for any $\gamma\in\Gamma$ and $x\in X_0$;
    \item[(2)] For any $x, x'\in X_0$ with $\Gamma\cdot x\neq\Gamma\cdot x'$ (i.e., $x$ and $x'$ belong to different $\Gamma$-orbits), we have $U_x\cap U_{x'}=\emptyset$;
    \item[(3)] For each $x\in X_0$, we have $U_x \subseteq O_{x,r}$, where $\{O_{x,r}\}$ is the coarse atlas family from Definition \ref{def:proper-Gamma-X-algebra}.
\end{enumerate}
\end{definition}

\vspace{-0.7cm}
\begin{figure}[htbp]
\centering
\small
\begin{tikzpicture}[
    scale=0.9,
    region/.style={dashed, thick},
    arrow/.style={->, >=stealth, thick},
    every node/.style={font=\footnotesize}
]

\fill[blue!3, rounded corners=3ex] (0,0) rectangle (6.5,5.8);
\node[font=\small\bfseries] at (3.25, 6.2) {Space $X$};

\draw[teal!60, fill=teal!8, thick, dashed, rounded corners=2ex] (2.3, 2.8) ellipse (1.8 and 1.15);
\node[teal!70] at (1.6, 1.4) {$\Gamma$-invariant countable subset: $X_0$};

\draw[orange!70, densely dashed, thick] (2.5, 3.1) -- (5.3, 4.2);
\node[orange!85] at (5.6, 4.4) {$\Gamma \cdot x_1$};
\fill[orange!85] (2.5, 3.1) circle (2.8pt) node[below=1pt, orange!90] {$x_1$};
\fill[orange!85] (3.3, 3.4) circle (2.5pt);
\fill[orange!85] (4.2, 3.8) circle (2.5pt);
\fill[orange!85] (5.0, 4.1) circle (2.5pt);

\draw[cyan!70, densely dashed, thick] (3.2, 2.3) -- (5.8, 1.5);
\node[cyan!85] at (6.1, 1.3) {$\Gamma \cdot x_2$};
\fill[cyan!85] (3.2, 2.3) circle (2.8pt) node[below=1pt, cyan!90] {$x_2$};
\fill[cyan!85] (4.0, 2.0) circle (2.5pt);
\fill[cyan!85] (4.8, 1.8) circle (2.5pt);
\fill[cyan!85] (5.6, 1.55) circle (2.5pt);

\draw[purple!70, densely dashed, thick] (2.7, 2.7) -- (5.0, 3.0);
\node[purple!85] at (5.2, 3.2) {$\Gamma \cdot x_3$};
\fill[purple!85] (2.7, 2.7) circle (2.8pt) node[below=1pt, purple!90] {$x_3$};
\fill[purple!85] (3.6, 2.9) circle (2.5pt);
\fill[purple!85] (4.3, 3.0) circle (2.5pt);

\draw[<->, red!80, thick] (5.2, 3.6) -- node[right=2pt, red!80, align=center] {different\\orbits} (5.5, 2.1);

\node at (8.2, 3.6) {\Large $\leadsto$};
\node[font=\itshape, align=center] at (8.2,3.2) {orbit separation\\condition};

\draw[->, gray!60, dashed, thick] (2.7, 3.2) to[bend left=15] (9.5, 4.6);
\draw[->, gray!60, dashed, thick] (3.3, 2.4) to[bend right=10] (9.5, 2.9);
\draw[->, gray!60, dashed, thick] (2.9, 2.8) to[bend right=5]  (9.5, 1.6);

\fill[cyan!3, rounded corners=3ex] (9.5,0) rectangle (17,5.8);
\node[font=\small\bfseries] at (13.35, 6.2) {Space $M$};

\draw[orange!50, dashed, fill=orange!3] (10.8, 4.6) ellipse (1.1 and 0.7);
\node[orange!80] at (10.2, 5.4) {$O_{x_1,r}$};
\draw[orange!50, dashed, fill=orange!3] (12.6, 4.6) ellipse (1.1 and 0.7);
\node[orange!80] at (12.0, 5.4) {$O_{\gamma x_1,r}$};
\draw[orange!50, dashed, fill=orange!3] (14.4, 4.6) ellipse (1.1 and 0.7);
\node[orange!80] at (13.8, 5.4) {$O_{\gamma' x_1,r}$};

\draw[orange!80, fill=orange!15, thick] (10.8, 4.6) ellipse (0.55 and 0.35);
\node[orange!90] at (10.6, 4.5) {$U_{x_1}$};
\draw[orange!80, fill=orange!15, thick] (12.6, 4.6) ellipse (0.55 and 0.35);
\node[orange!90] at (12.5, 4.5) {$U_{\gamma x_1}$};
\draw[orange!80, fill=orange!15, thick] (14.4, 4.6) ellipse (0.55 and 0.35);
\node[orange!90] at (14.2, 4.5) {$U_{\gamma' x_1}$};

\draw[->, >=stealth, red!60, thick] (11.4, 4.7) -- node[above=-2pt, scale=0.8] {$\gamma$} (12.0, 4.7);
\draw[->, >=stealth, red!60, thick] (13.2, 4.7) -- node[above=-2pt, scale=0.8] {$\gamma'$} (13.8, 4.7);

\draw[cyan!50, dashed, fill=cyan!3] (11.6, 3.35) ellipse (1.1 and 0.7);
\node[cyan!80] at (11.6, 4.0) {$O_{x_2,r}$};
\draw[cyan!50, dashed, fill=cyan!3] (13.4, 3.35) ellipse (1.1 and 0.7);
\node[cyan!80] at (13.4, 4.0) {$O_{\gamma x_2,r}$};
\draw[cyan!50, dashed, fill=cyan!3] (15.2, 3.35) ellipse (1.1 and 0.7);
\node[cyan!80] at (15.2, 4.0) {$O_{\gamma' x_2,r}$};

\draw[cyan!80, fill=cyan!15, thick] (11.6, 3.35) ellipse (0.55 and 0.35);
\node[cyan!90] at (11.6, 3.2) {$U_{x_2}$};
\draw[cyan!80, fill=cyan!15, thick] (13.4, 3.35) ellipse (0.55 and 0.35);
\node[cyan!90] at (13.4, 3.2) {$U_{\gamma x_2}$};
\draw[cyan!80, fill=cyan!15, thick] (15.2, 3.35) ellipse (0.55 and 0.35);
\node[cyan!90] at (15.2, 3.2) {$U_{\gamma' x_2}$};

\draw[->, >=stealth,red!60, thick] (12.2, 3.3) -- node[above=-2pt, scale=0.8] {$\gamma$} (12.8, 3.3);
\draw[->, >=stealth,red!60, thick] (14.0, 3.3) -- node[above=-2pt, scale=0.8] {$\gamma'$} (14.6, 3.3);

\draw[purple!50, dashed, fill=purple!3] (10.8, 1.6) ellipse (1.1 and 0.7);
\node[purple!80] at (10.2, 2.4) {$O_{x_3,r}$};
\draw[purple!50, dashed, fill=purple!3] (12.6, 1.6) ellipse (1.1 and 0.7);
\node[purple!80] at (12.0, 2.4) {$O_{\gamma x_3,r}$};
\draw[purple!50, dashed, fill=purple!3] (14.4, 1.6) ellipse (1.1 and 0.7);
\node[purple!80] at (13.8, 2.4) {$O_{\gamma' x_3,r}$};

\draw[purple!80, fill=purple!15, thick] (11.1, 1.55) ellipse (0.8 and 0.35);
\node[purple!90] at (10.7, 1.15) {$U_{x_3}$};
\draw[purple!80, fill=purple!15, thick] (12.6, 1.55) ellipse (0.8 and 0.35);
\node[purple!90] at (12.6, 1.15) {$U_{\gamma x_3}$};
\draw[purple!80, fill=purple!15, thick] (14.1, 1.55) ellipse (0.8 and 0.35);
\node[purple!90] at (14.5, 1.15) {$U_{\gamma' x_3}$};

\draw[->, >=stealth, red!60, thick] (11.4, 1.55) -- node[above=-2pt, scale=0.8] {$\gamma$} (12.0, 1.55);
\draw[->, >=stealth, red!60, thick] (13.2, 1.55) -- node[above=-2pt, scale=0.8] {$\gamma'$} (13.8, 1.55);

\node[red!80, scale=0.95, font=\footnotesize] at (12.6, 0.82) {$U$'s within the same orbit may overlap};

\draw[<-, red!70, thick] (12.1, 4.0) to[bend right=15] (12.6, 5.0) node[right, red!70, font=\footnotesize] {$O$'s might overlap};

\draw[<->, red!90, thick] (14.4, 4.2) -- (14.4, 3.75);
\node[red!90, right, align=left] at (14.5, 4.0) {orbits  are strictly\\disjoint};

\end{tikzpicture}

\caption{The family $\{U_x\}_{x \in X_0}$ separates the distinct $\Gamma$-orbits of $X_0$.}
\label{fig:separated-family}
\end{figure}
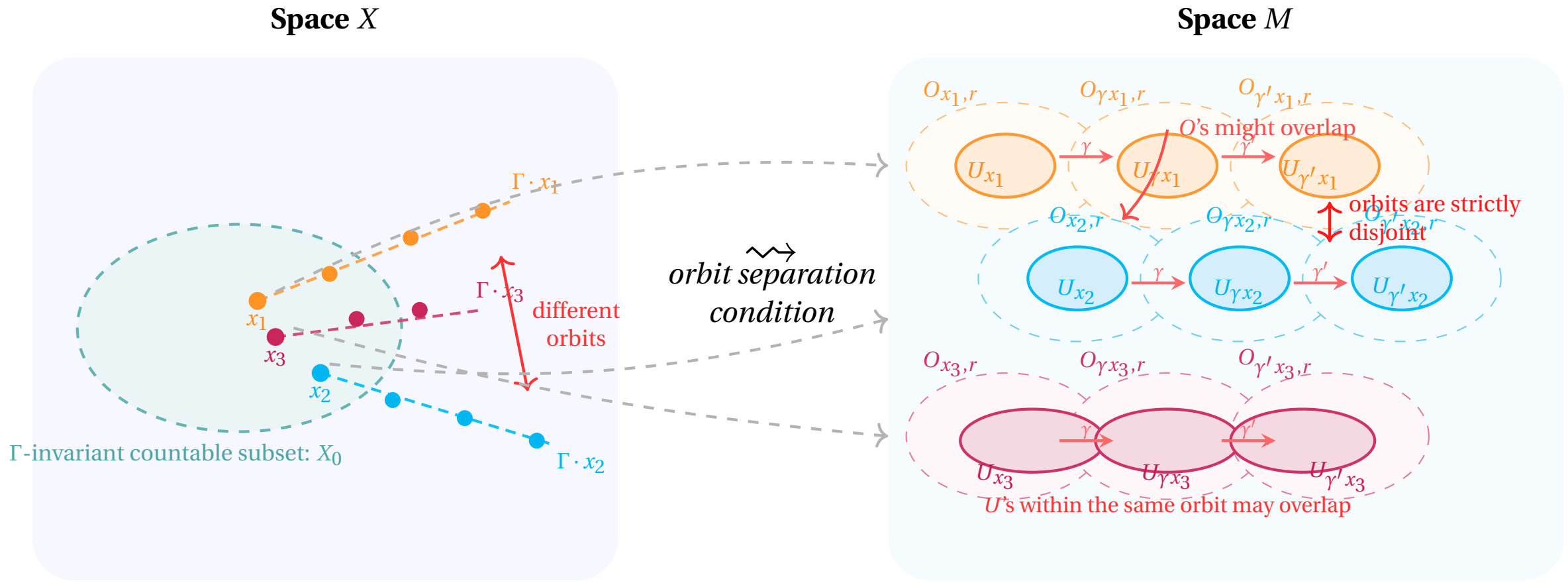

\vspace{-0.8cm}
\begin{remark}
Definition \ref{definition 3.10} provides a way to decompose the coefficient space $M$ into localized pieces that respect both the $\Gamma$-action and the coarse geometry of $X$. Condition (1) ensures the decomposition is $\Gamma$-equivariant. Condition (2) ensures that pieces belonging to different $\Gamma$-orbits are disjoint, which allows us to localize operators orbit-by-orbit. Condition (3) places each piece $U_x$ within the coarse atlas $O_{x,r}$, thereby linking the topological decomposition of $M$ back to the proper $(\Gamma,X)$-structure of $\mathcal{A}$. Together, these conditions allow us to analyze different orbits separately, while keeping everything localized and compatible with the coarse geometric data of the underlying space.
See Fig. \ref{fig:separated-family} for an illustration of this construction. However, the open subsets indexed by a single orbit are allowed to overlap. Under the additional conditions in our next example, one obtains a refinement with pairwise disjoint translates. We shall consider this case first (Proposition \ref{proposition:eval-slice-separated}) and then generalize it to an arbitrary $(\Gamma, X_0,r)$-separated family (Theorem \ref{thm:evaluation-map-isomorphism}).
\end{remark}

\begin{example}\label{example:slice-separated-sets}
Let $r>0$ and $X_0\subseteq X$ be a $\Gamma$-invariant
countable subset. Choose a set of orbit representatives
$\{x_i\}_{i\in I}$ for $X_0/\Gamma$. For each $i\in I$, let
$$F_i:=\Gamma_{x_i}
=\{\gamma\in\Gamma:\gamma x_i=x_i\}$$
denote the stabilizer subgroup of $x_i$.
For each $i\in I$, suppose that there is an open subset
$U_{i,0}\subseteq M$ satisfying the following conditions:
\begin{enumerate}
\item[(1)] $U_{i,0}\subseteq O_{x_i,r}$ and $U_{i,0}$ is
$F_i$-invariant, i.e.,
$gU_{i,0}=U_{i,0}$ for all $g\in F_i$.
\item[(2)] $U_i:=\Gamma \cdot U_{i,0}$ are pairwise disjoint, i.e., $U_i \cap U_j=\varnothing$ whenever $i\neq j$.
\item[(3)] For every $\gamma,\gamma'\in\Gamma$, if $\gamma U_{i,0}\cap \gamma' U_{i,0}\neq\varnothing$, then $\gamma F_i=\gamma'F_i$. Equivalently, if $\gamma U_{i,0}\cap U_{i,0}\neq\varnothing$, then $\gamma\in F_i$.
\end{enumerate}
Condition (3), together
with the $F_i$-invariance of $U_{i,0}$, gives the disjoint union
\[U_i=\bigsqcup_{[\gamma]\in\Gamma/F_i}\gamma U_{i,0}.\]
Consequently, the map
$\Gamma\times_{F_i}U_{i,0}
\rightarrow U_i$,
$[\gamma,u]\mapsto\gamma u$,
is a well-defined $\Gamma$-equivariant homeomorphism. Here
$\Gamma\times_{F_i}U_{i,0}$ is the quotient of
$\Gamma\times U_{i,0}$ by the relation
$(\gamma g,u)\sim(\gamma,gu)$,
$\gamma\in\Gamma$, $g\in F_i$, $u\in U_{i,0}$.

For every $x\in X_0$, there exist $i\in I$ and 
$\gamma\in\Gamma$ such that
$x=\gamma x_i$.
The index $i$ is unique, while $\gamma$ is determined only up to right
multiplication by an element of $F_i$. We define
\[U_x:=\gamma U_{i,0}.\]
The family
$\{U_x\}_{x\in X_0}$ is $(\Gamma,X_0,r)$-slice separated in the sense
of Definition \ref{definition 3.10}.
Indeed, if $x=\gamma x_i$, then for every $\eta\in\Gamma$, $$U_{\eta x} =U_{\eta\gamma x_i}=\eta\gamma U_{i,0}=\eta U_x.$$ Thus, the family $\{U_x\}_{x\in X_0}$ is $\Gamma$-equivariant. Since $U_{i,0}\subseteq O_{x_i,r}$ and the collection $\{O_{x,r}\}$ is $\Gamma$-equivariant, $$U_x=\gamma U_{i,0}\subseteq\gamma O_{x_i,r}=O_{\gamma x_i,r}=O_{x,r}.$$ Next, we
suppose that $x=\gamma x_i$ and $x'=\gamma' x_j$ belong to different $\Gamma$-orbits. Then $i\neq j$, and hence $U_x\subseteq U_i$ and $U_{x'}\subseteq U_j$ are disjoint. This completes the verification. 
Moreover, if $x'=\eta x$ and $x'\neq x$, then $U_x\cap U_{x'}=\varnothing$. setting $x=\gamma x_i$, we have $x'=\eta\gamma x_i$, and therefore
$U_x=\gamma U_{i,0}$ and $U_{x'}=\eta\gamma U_{i,0}$. If these two sets
intersected, condition (3) would imply $\eta\gamma F_i=\gamma F_i$, hence
$\gamma^{-1}\eta\gamma\in F_i$. This would give
$\eta x=\eta\gamma x_i=\gamma(\gamma^{-1}\eta\gamma)x_i=\gamma x_i=x$,
contradicting $x'\neq x$.

We call the collection 
$\{U_i=\Gamma\cdot U_{i,0}\}_{i\in I}$
a \emph{$\boldsymbol{\Gamma}$\textbf{-slice separated family}}, and 
$U_{i,0}$ the \emph{\textbf{local slices}}.
\end{example}
The $\Gamma$-slice separated families are tailor-made for the geometric cutting-and-pasting technique that will be central to the proof of Theorem \ref{main theorem}. By localizing the problem to $\Gamma$-slices with finite symmetries, we can reduce the original equivariant $K$-theory computation on the infinite group $\Gamma$ to a family of finite group problems. To incorporate the geometry of $X$ into this decomposition, we now geometrically pair each local slice $U_{i,0}$ with a closed subset $Y_i$ in the Rips complex $P_{d}(X)$ through the common finite subgroup $F_i$ and the representative point $x_i$.

Let $\{x_i\}_{i\in I}$ be the chosen complete set of orbit representatives in $X_0$ for an index set $I$. For each $i \in I$, let $F_i=\Gamma_{x_i}$ be its finite stabilizer subgroup. Fix $d>0$ and consider a family $\{Y_i\}_{i\in I}$ of closed subsets of $P_{d}(X)$ satisfying:
\begin{enumerate}
\item[(1)] $x_i\in Y_i$ for each $i\in I$;
\item[(2)] $Y_i$ is $F_i$-invariant for each $i\in I$;
\item[(3)] there exists a uniform radius $r_0>0$, independent of $i$, such that $Y_i\subseteq \overline{B}_{P_{d}(X)}(x_i,r_0)$ for all $i\in I$.
\end{enumerate}
Three such families are:
\begin{enumerate}
\item[(a)] $Y_i=\{x_i\}$ for each $i\in I$.
\item[(b)] $Y_i=\overline{B}_{P_{d}(X)}(x_i,R)$, the closed ball of radius $R > 0$ centered at $x_i$;
\item[(c)] If $d \geq 2R$, 
let $Y_i=\Delta_i^{(d)}(R)$ be the simplex with vertex set $$V_i(R):=\{v\in X\mid d_X(v,x_i)\leq R\}.$$ 
\end{enumerate}
Note that for any $v,w\in V_i(R)$, $d_X(v,w)\leq d_X(v, x_i)+d_X(x_i,w)\leq 2R\leq d$, so $V_i(R)$ spans a simplex in $P_{d}(X)$. Since $F_i$ fixes $x_i$ and acts isometrically on $X$, $\Delta_i^{(d)}(R)$ is $F_i$-invariant. Since the simplex in $P_{d}(X)$ is bounded by its vertex diameter, we have $\Delta_i^{(d)}(R)\subseteq \overline{B}_{P_{d}(X)}(x_i,r_0)$ for a uniform $r_0>0$.

Now, let $\{U_i\}_{i\in I}$ be a $\Gamma$-slice separated family as in Example \ref{example:slice-separated-sets}, where $U_i=\Gamma\cdot U_{i,0}$. Recall that $\phi:C_0(M)\to Z(M(\mathcal{A}))$ denotes the $\Gamma$-equivariant $*$-homomorphism from Definition \ref{def:proper-Gamma-X-algebra}(1). For an open subset $U\subseteq M$ we regard $C_0(U)$ as an ideal of $C_0(M)$, and set 
$$\mathcal{A}_U:=\overline{\operatorname{span}}\{\phi(f)a\mid f\in C_0(U), a \in\mathcal{A}\}.$$ 
Then $\mathcal{A}_U$ is a closed two sided ideal of $\mathcal{A}$. The equivariance of $\phi$ implies $\gamma(\phi(f)a)=\phi(\gamma f)(\gamma a)$ for all $\gamma\in\Gamma$, whence $\gamma\cdot\mathcal{A}_U=\mathcal{A}_{\gamma U}$. We apply this construction to $U_i$ and $U_{i,0}$. Since $U_i$ is $\Gamma$-invariant, $\mathcal{A}_{U_i}$ is a $\Gamma$-invariant ideal of $\mathcal{A}$. Since $U_{i,0}$ is $F_i$-invariant by condition (1) of Example \ref{example:slice-separated-sets}, $\mathcal{A}_{U_{i,0}}$ is invariant under the restricted action of $F_i$ and is thus an $F_i$-$C^*$-algebra. 

Let $\{Y_i\}_{i\in I}$ be a family of closed $F_i$-invariant subsets of $P_{d}(X)$ as above. Since $\mathcal{A}_{U_{i,0}}$ is an $F_i$-invariant ideal of $\mathcal{A}$ and $\phi$ restricts to $C_0(U_{i,0}) \to Z(M(\mathcal{A}_{U_{i,0}}))$, the sets $$O^{(i)}_{z,r}:=O_{z,r}\cap U_{i,0}$$ endow $\mathcal{A}_{U_{i,0}}$ with a proper $(F_i,X)$-structure inheriting the same control functions $\rho_\pm$. By restricting the $\Gamma$-action to $F_i$ in Definition \ref{def:equiv-alg-roe-X}, we obtain the $F_i$-equivariant Roe algebra $C^*(Y_i,\mathcal{A}_{U_{i,0}})^{F_i}$ and its localization algebra $C_L^*(Y_i,\mathcal{A}_{U_{i,0}})^{F_i}$. We still use the module $\ell^2(\Gamma)\otimes H_0$ with the $F_i$-action. By equivariant bounded geometry (Definition \ref{def:equiv-bdd-geom}, Remark \ref{rem:equivariant-bounded-geometry}), $\text{sup}_{i\in I}|F_{i}|\leq N<\infty$, the family $\{C^*(Y_i,\mathcal{A}_{U_{i,0}})^{F_i}\}_{i\in I}$ is uniformly controlled.
\begin{definition}\label{def:product-algebra}
Let $\{U_i\}_{i\in I}$ be a $\Gamma$-slice separated family as in Example \ref{example:slice-separated-sets}, with $U_i=\Gamma\cdot U_{i,0}$. Fix a $\Gamma$-invariant dense subset $Z_{d}\subseteq P_{d}(X)$. For $i\in I$, set $Z_i:=Z_{d}\cap Y_i$, which is assumed to be dense in $Y_i$. 
\begin{enumerate}
\item[(1)]
The \textit{uniformly controlled algebraic product Roe algebra}, denoted by $\mathbb A_{\mathrm u}[Y_i,i\in I]$, is the subalgebra of $\prod_{i\in I}\mathbb C[Y_i,\mathcal{A}_{U_{i,0}}]^{F_i}$ consisting of elements $(T_i)_{i\in I}$ such that
$$\sup_{i\in I}\|T_i\|<\infty,\ \sup_{i\in I}\operatorname{prop}(T_i)<\infty,\ \sup_{i\in I}\operatorname{row}(T_i)<\infty,\ \sup_{i\in I}\operatorname{col}(T_i)<\infty.$$
The product Roe algebra associated with $\{(Y_i,U_{i,0},F_i)\}_{i\in I}$ is the norm closure\[
\mathcal{A}^*(Y_i,i\in I):=
\overline{\mathbb{A}_{\mathrm u}[Y_i,i\in I]}^{\|\cdot\|_\infty}
\subseteq
\prod_{i\in I}
C^*(Y_i,\mathcal{A}_{U_{i,0}})^{F_i}.
\]
\item[(2)]
Let $\mathbb{A}_L^*(Y_i,i\in I)$ be the $*$-algebra of all
bounded, uniformly norm continuous functions
\[b:[0,\infty)\longrightarrow
\mathcal{A}^*(Y_i,i\in I),
\quad
b(t)=(b_i(t))_{i\in I},\]
satisfying
\begin{enumerate}
\item[(i)]
$b(t)\in\mathbb{A}_{\mathrm u}[Y_i,i\in I]$ for every $t\geq 0$;
\item[(ii)]the propagation is uniformly bounded, i.e.,
$\sup_{t\geq0}\sup_{i\in I}
\operatorname{prop}(b_i(t))<\infty$;
\item[(iii)]the propagation uniformly tends to zero, i.e.,
$\lim_{t\to\infty}
\sup_{i\in I}\operatorname{prop}(b_i(t))=0$.
\end{enumerate}
The product localization algebra is defined as
\[\mathcal{A}_L^*(Y_i,i\in I):=
\overline{\mathbb A_L^*(Y_i,i\in I)}
\subseteq
C_{b,u}\bigl([0,\infty),
\mathcal{A}^*(Y_i,i\in I)\bigr),\]
where
\[\|b\|=\sup_{t\geq0}\|b(t)\|=\sup_{t\geq0,\ i\in I}\|b_i(t)\|.
\]
\end{enumerate}
\end{definition}
Then there is an evaluation-at-zero map:
\[e:\mathcal{A}_L^{*}(Y_i,i\in I)\longrightarrow \mathcal{A}^{*}(Y_i,i\in I),\quad e(b)=b(0)=(b_{i}(0))_{i\in I}.\]

\begin{lemma}\label{lemma 3.16}
Let $\{U_i\}_{i\in I}$ be a $\Gamma$-slice separated family
as in Example \ref{example:slice-separated-sets}, and put
$U:=\bigcup_{i\in I}U_i$.
For each $d>0$, let
$C^{*}(P_{d}(X),\mathcal{A})^\Gamma_U$
and
$C_L^{*}(P_{d}(X),\mathcal{A})^\Gamma_U$
denote the closed ideals obtained by restricting the coefficients
to
$\mathcal{A}_U:=\overline{C_0(U)\mathcal{A}}$.
Set
$\mathcal{I}:=
\{(d,R)\in(0,\infty)^2\mid d\geq2R\}$. Then there are canonical 
isomorphisms:
\begin{enumerate}
\item[(1)] $C^{*}(P_{d}(X),\mathcal{A})^{\Gamma}_{U} \cong \lim\limits_{R\to\infty} \mathcal{A}^{*}\bigl(\overline{B}_{P_{d}(X)}(x_{i},R),\; i\in I\bigr)$;
\item[(2)] $C^{*}_{L}(P_{d}(X),\mathcal{A})^{\Gamma}_{U}\cong \lim\limits_{R\to\infty} \mathcal{A}^{*}_{L}\bigl(\overline{B}_{P_{d}(X)}(x_{i},R),\ i\in I\bigr)$;
\item[(3)] $\lim\limits_{d\to\infty} C^{*}(P_{d}(X),\mathcal{A})^{\Gamma}_{U} \cong \lim\limits_{(d,R)\in \mathcal{I}} \mathcal{A}^{*}\bigl(\Delta_i^{(d)}(R),\ i\in I\bigr)$;
\item[(4)] $\lim\limits_{d\to\infty}C_{L}^{*}(P_{d}(X),\mathcal{A})^{\Gamma}_{U} \cong \lim\limits_{(d,R)\in \mathcal{I}} \mathcal{A}^{*}_{L}\bigl(\Delta_i^{(d)}(R),\ i\in I\bigr)$.
\end{enumerate}
\end{lemma}
\begin{proof}[\textbf{Proof.}]\normalsize
We first prove (1). For open subset $U$ and a fixed scale $d>0$, let $Z_{d}$ be a $\Gamma$-invariant countable dense subset of the Rips complex $P_{d}(X)$. We define the Hilbert $\mathcal{A}$-module $H_{P_{d}(X),\mathcal{A},U}$ by$$H_{P_{d}(X),\mathcal{A},U}:=\ell^{2}(Z_{d}) \otimes \ell^{2}(\Gamma)\otimes H_{0}\otimes\mathcal{A}_{U}.$$
It is a submodule of $H_{P_{d}(X),\mathcal{A}}$, and the closed ideal $C^{*}(P_{d}(X),\mathcal{A})_U^\Gamma$ has a faithful non-degenerate representation on $H_{P_{d}(X),\mathcal{A},U}$.
Since $\{U_i\}_{i\in I}$ is a family of disjoint $\Gamma$-slice separated sets with $U_i\cong\Gamma \times_{F_i} U_{i,0}$, we have the following direct sum decomposition:
\[\mathcal{A}_U=\bigoplus_{i\in I}\mathcal{A}_{U_i}=\bigoplus_{i\in I}\bigoplus_{[\gamma F_{i}]\in\Gamma/F_i}\mathcal{A}_{\gamma\cdot U_{i,0}}.\] 
The Hilbert $\mathcal{A}$-module splits accordingly:
\[H_{P_{d}(X),\mathcal{A},U}=\bigoplus_{i\in I}\bigoplus_{[\gamma F_{i}]\in\Gamma/F_i} H_{P_{d}(X),\mathcal{A},\gamma\cdot U_{i,0}}.\]
For any $T\in\mathbb{C}[P_{d}(X),\mathcal{A}]_U^\Gamma$, this splitting decomposes the operator into
\[T=\bigoplus_{i\in I}\bigoplus_{[\gamma F_{i}]\in \Gamma/F_i}T_{\gamma,i},\]
where \[T_{\gamma,i}(x,y):=(\mathrm{id}\otimes p_{\gamma,i})\bigl(T(x,y)\bigr)\in\mathcal{K}(\ell^2(\Gamma)\otimes H_0)\otimes\mathcal{A}_{\gamma\cdot U_{i,0}},\]
with $p_{\gamma,i}:\mathcal{A}_U\to \mathcal{A}_{\gamma\cdot U_{i,0}}$ being the projection induced by the $c_0$-direct sum decomposition.
For any $\gamma,\gamma'\in\Gamma$, we have that
\[T_{\gamma,i}(\gamma'x,\gamma'y)=(\mathrm{id}\otimes p_{\gamma,i})\bigl(T(\gamma' x,\gamma' y)\bigr)=\gamma'\cdot(\mathrm{id}\otimes p_{\gamma'^{-1}\gamma,i})\bigl(T(x, y)\bigr)=\gamma'\cdot T_{\gamma'^{-1}\gamma,i}(x, y).\]
By substituting $\gamma'=\gamma$, we get
$T_{\gamma,i}(\gamma x,\gamma y)=\gamma\cdot T_{e,i}(x, y)$,
where $e$ denotes the identity element of $\Gamma$. Then we have that
\[T=\bigoplus_{i\in I}\bigoplus_{[\gamma F_{i}]\in \Gamma/F_i}T_{\gamma,i}=\bigoplus_{i\in I}\bigoplus_{[\gamma F_{i}]\in\Gamma/F_i}\gamma T_{e,i}.\]

Since $T$ is $\Gamma$-equivariant and $U_{i,0}$ is $F_i$-invariant, for any $\gamma\in F_i$ we have
\[T_{e,i}(\gamma x,\gamma y)=(\mathrm{id}\otimes p_{e,i})\bigl(T(\gamma x,\gamma y)\bigr)=\gamma\cdot (\mathrm{id}\otimes p_{e,i})\bigl(T(x,y)\bigr)=\gamma\cdot T_{e,i}(x,y).\]
This implies that each $T_{e,i}$ is $F_i$-equivariant. By condition (1) of Example \ref{example:slice-separated-sets}, 
$\operatorname{Supp}(T_{e,i}(x, y))\subseteq U_{i,0}\subseteq O_{x_i,r}$ for a fixed $r > 0$.
Moreover, by Definition \ref{def:equiv-alg-roe-X}(5), there exists $r'(T)>0$, such that
$\operatorname{Supp}(T(x, y))\subseteq O_{x,r'(T)}$.
Recall that for $x\in Z_{d}\subseteq P_{d}(X)$, the set $O_{x,r'(T)}$ is defined as $\bigcup_{v\in V_x} O_{v,r'(T)}$, where $V_x\subset X$ is the set of vertices of the simplex containing $x$. Choose $s>0$ large enough so that
$\rho_-(s)\geq\max\{r,r'(T)\}$.
If $T_{e,i}(x, y)\neq 0$, then
$O_{x,r'(T)}\cap O_{x_i,r}\neq\emptyset$.
This implies that there exists some $v_x\in V_x$ such that
$O_{v_x,r'(T)}\cap O_{x_i,r}\neq\emptyset$.
By inclusion,
$O_{v_x,\rho_-(s)}\cap O_{x_i,\rho_-(s)}\neq\emptyset$.
Applying Definition \ref{def:proper-Gamma-X-algebra}(2)(b), we obtain
$d_X(v_x, x_i)\leq 2s$. Since the natural inclusion $X\hookrightarrow P_{d}(X)$ is bornologous, there exists a constant $C_{d}>0$ such that
\[ d_{P_{d}(X)}(v_x, x_i)\leq C_{d}\cdot d_X(v_x, x_i)+C_{d}\leq 2sC_{d}+C_{d}.\]
Since $X$ has bounded geometry, the simplices of $P_{d}(X)$ have uniformly bounded diameter. Thus, for $v_x$ a vertex of the simplex containing $x$ there is a constant $\kappa=\kappa(X,d)$ with $d_{P_{d}(X)}(x,v_x)\leq\kappa$. If $T$ has finite propagation $L$, any non-zero entry requires $d_{P_{d}(X)}(x,y)\leq L$. We obtain \[d_{P_{d}(X)}(y, x_i)\leq d_{P_{d}(X)}(y,x)+d_{P_{d}(X)}(x,v_x)+d_{P_{d}(X)}(v_x,x_i)\leq L+\kappa+2sC_{d}+C_{d}.\] Similarly, $d_{P_{d}(X)}(x,x_i)\leq\kappa+2sC_{d}+C_{d}$. Therefore, setting $R:=R(T,d):=L+\kappa+2sC_{d}+C_{d}$, we conclude \[T_{e,i}\in\mathbb{C}[\overline{B}_{P_{d}(X)}(x_i,R),\mathcal{A}_{U_{i,0}}]^{F_i}.\]

Denote by $N_{\mathrm r}(T)$ and $N_{\mathrm c}(T)$ the row and column
bounds for $T$ as in Definition
\ref{def:equiv-alg-roe-X}(2). Then
$\operatorname{row}(T_{e,i})\leq N_{\mathrm r}(T)$,
$\operatorname{col}(T_{e,i})\leq N_{\mathrm c}(T)$
for every $i\in I$. Moreover,
$\operatorname{prop}(T_{e,i})
\leq \operatorname{prop}(T)$.
Hence,
\[(T_{e,i})_{i\in I}\in
\mathbb{A}_{\mathrm u}
\bigl[\overline B_{P_{d}(X)}(x_i,R),i\in I\bigr].\]
In addition, note that \begin{equation*}
\|T\|=\sup_{i\in I,[\gamma]\in\Gamma/F_i}\|T_{\gamma,i}\|=\sup_{i\in I}\|T_{e,i}\|<\infty.
\end{equation*}
Consequently,
the map $T\mapsto (T_{e,i})_{i\in I}$ defines an isometric $*$-homomorphism from $\mathbb C[P_{d}(X),\mathcal{A}]^\Gamma_U$ into the direct limit $\lim_{R\to\infty}\mathcal{A}^*\bigl(\overline B_{P_{d}(X)}(x_i,R), i\in I\bigr)$. 
Therefore, it extends uniquely to an injective $C^*$-homomorphism upon completion, given by
\[\Phi: C^*(P_{d}(X),\mathcal{A})^\Gamma_U \longrightarrow\lim_{R\to\infty}\mathcal{A}^*\bigl( \overline B_{P_{d}(X)}(x_i,R), i\in I \bigr),\quad \Phi(T)=(T_{e,i})_{i\in I}.\]\par
Conversely, let
$(T_i)_{i\in I}\in
\mathbb{A}_{\mathrm u}
\bigl[\overline B_{P_{d}(X)}(x_i,R),i\in I\bigr].
$
Choose $N>0$ such that
$\operatorname{row}(T_i)\leq N$,
$\operatorname{col}(T_i)\leq N$
for all $i\in I$. For $i\in I$ and $[\gamma]\in\Gamma/F_i$, we define
\[T_{\gamma,i}(x,y)=
\gamma\cdot
T_i(\gamma^{-1}x,\gamma^{-1}y),
\quad x,y\in Z_{d}.\]
The $F_i$-equivariance of $T_i$ ensures that this definition is independent of the choice of representative for $[\gamma]$. We then define
\[
\Psi\bigl((T_i)_{i\in I}\bigr):=
\bigoplus_{i\in I}
\bigoplus_{[\gamma]\in\Gamma/F_i}
T_{\gamma,i}.\]
By the direct sum decomposition of $H_{P_{d}(X),\mathcal{A},U}$ and $\sup_{i\in I}\Vert{}T_i\Vert{}<\infty$, $\Psi\bigl((T_i)_{i\in I}\bigr)$ is a well-defined bounded block-diagonal operator on $H_{P_{d}(X),\mathcal{A},U}$. It is $\Gamma$-equivariant and has finite propagation at most $2R$. Moreover, if $\gamma^{-1}x \in \overline B_{P_{d}(X)}(x_i,R)$, the coarse equivalence $X\hookrightarrow P_{d}(X)$ yields a non-decreasing function $q_{d}$ such that every vertex $v$ of the simplex containing $\gamma^{-1}x$ satisfies 
\[d_X(v,x_i)\leq q_{d}\bigl(d_{P_{d}(X)}(v,x_{i})\bigr)\leq q_{d}\bigl(d_{P_{d}(X)}(v,\gamma^{-1}x)+d_{P_{d}(X)}(\gamma^{-1}x, x_{i})\bigr)\leq q_{d}(\kappa+R),\]
where $\kappa=\kappa(X,d)$ is the uniform bound for the distance from any point in a simplex of $P_{d}(X)$ to its vertices. Furthermore, 
since $U_{i,0}\subseteq O_{x_i,r}$ for some $r>0$, we have that
\[U_{i,0}\subseteq O_{x_i,r}
\subseteq O_{v,\ \rho_{+}(d_X(v,x_i))+r}\subseteq O_{v,\ r_R},
\ \text{where}\
r_R:=r+\rho_{+}\bigl(q_{d}(\kappa+ R)\bigr).
\]
This implies $U_{i,0}\subseteq O_{\gamma^{-1}x, r_R}$. We obtain
$\gamma\cdot U_{i,0}\subseteq\gamma\cdot O_{\gamma^{-1}x,r_R}=O_{x, r_R}$.

Since
$Y_i\subseteq\overline B_{P_{d}(X)}(x_i,R)$,
there exists a constant $M_{d,R}>0$ such that every $x\in Z_{d}$
belongs to at most $M_{d,R}$ of the translated subsets
$\gamma Y_i$, where $i\in I$ and $[\gamma]\in\Gamma/F_i$.
Indeed, if $x\in\gamma Y_i$, then
$d_{P_{d}(X)}(x,\gamma x_i)\leq R$.
Let $v_x$ be a vertex of the simplex containing $x$. Then we have that
\[d_{P_{d}(X)}(v_x,\gamma x_i)\leq\kappa+R.\]
So
\[d_X(v_x,\gamma x_i)\leq q_{d}(\kappa+R).\]
Since $X$ has bounded geometry,
\[M_{d,R}\leq\sup_{v\in X}\#B_X\bigl(v,q_{d}(\kappa+R)\bigr)
<\infty.
\]
Then we have that
\[
\text{row}\bigl(\Psi((T_i))\bigr)
\leq
M_{d,R}\sup_{i}\text{row}(T_i)<\infty,
\]
and similarly,
\[
\text{col}\bigl(\Psi((T_i))\bigr)
\leq
M_{d,R}\sup_i\text{col}(T_i)<\infty.
\]
Therefore, 
\[
\Psi\bigl((T_i)_{i\in I}\bigr)\in\mathbb C[P_{d}(X),\mathcal{A}]^\Gamma_U.
\]
It is clear that $\Psi$ is isometric. It extends continuously to the $C^{*}$-completion, and we obtain
\[\Psi:\lim_{R\to\infty}\mathcal{A}^{*}\bigl(\overline{B}_{P_{d}(X)}(x_i,R),i\in I\bigr)\longrightarrow C^{*}(P_{d}(X),\mathcal{A})_{U}^\Gamma.\]
A direct verification on the dense finite-propagation subalgebras confirms that $\Phi$ and $\Psi$ are inverses. Being isometric, they extend to $C^{*}$-isomorphisms, proving (1). Statement (2) is similar.

Applying (1) and (2) and then passing to the inductive limit as $d\to\infty$, we only need to compare $\overline{B}_{P_{d}(X)}(x_i,S)$ with $\Delta_i^{(d)}(R)$, where $S$ denotes the radius parameter appearing in (1). Let $D_{d}$ be a uniform upper bound for the diameters of the simplices in $P_{d}(X)$. If $(d,R)\in\mathcal I$, then the vertices contained in $\overline{B}_X(x_i,R)$ have pairwise distance at most $2R\leq d$, and hence span the simplex $\Delta_i^{(d)}(R)$. Since $x_i$ is one of these vertices, we obtain $\Delta_i^{(d)}(R)\subseteq\overline{B}_{P_{d}(X)}(x_i,D_{d})$. Conversely, for any $d,S>0$, the coarse equivalence between $P_{d}(X)$ and $X$ yields a constant $q_{d}(S)>0$, independent of $i$, such that all vertices of the smallest simplex containing any point of $\overline{B}_{P_{d}(X)}(x_i,S)$ lie in $\overline{B}_X(x_i,q_{d}(S))$. Thus, for any $d'\geq\max\{d,2q_{d}(S)\}$, the inclusion $P_{d}(X)\subseteq P_{d'}(X)$ gives $\overline{B}_{P_{d}(X)}(x_i,S)\subseteq\Delta_i^{(d')}(q_{d}(S))$. As these inclusions are uniform in $i$ and compatible with the structure maps, the two inductive systems have canonically isomorphic limits, proving (3) and (4).
\end{proof}

\section{Slicewise Evaluation Isomorphisms}

In \cite{Yu97}, Yu proved that the $K$-theory of localization algebras is invariant under strong Lipschitz homotopy equivalence. We now present an equivariant version of this result suitable for our context. In what follows, $I$ denotes a countable index set.
\begin{definition}\label{def:equivariant-lipschitz}
Let $\{X_i\}_{i\in I}$ and $\{Y_i\}_{i\in I}$ be two families of uniformly bounded closed subspaces of $P_{d}(X)$ at scale $d>0$. Suppose that for each $i\in I$, $X_i$ and $Y_i$ are invariant under a finite subgroup $F_i\leq\Gamma$, and have a common basepoint $x_i\in X_{i}\cap Y_i$.

A map $g:\bigsqcup_{i\in I}X_i\to\bigsqcup_{i\in I}Y_i$ satisfying $g(X_i)\subseteq Y_i$ for each $i\in I$ is called \emph{$F_i$-equivariant} if for each $i \in I$, all $x\in X_i$ and all $\gamma\in F_i$, one has $g(\gamma\cdot x)=\gamma\cdot g(x)$.
Furthermore, such a map $g$ is called \emph{uniformly Lipschitz} if 
there exists a constant $C>0$, independent of $i\in I$, such that $d(g(x),g(y))\leq Cd(x,y)$ for all $i \in I$ and all $x,y\in X_i$.
\end{definition}
\begin{definition}\label{def:equivariant-strong-lipschitz-homotopy}
Let $g_0,g_1:\bigsqcup_{i\in I}X_i\to\bigsqcup_{i\in I}Y_i$ be two $F_i$-equivariant uniformly Lipschitz maps. 
We say $g_0$ and $g_1$ are \emph{$F_i$-equivariant strongly Lipschitz homotopic} if there exists a map
\[H:\bigsqcup_{i\in I} X_i\times[0, 1]\to\bigsqcup_{i\in I}Y_i\]
satisfying:
\begin{enumerate}
\item [(1)] $H(x,0)=g_0(x)$ and $H(x,1)=g_1(x)$ for all $x\in\bigsqcup_{i\in I} X_i$;
\item [(2)] For each $t\in [0, 1]$, the map $H_t:= H(\cdot,t)$ is $F_i$-equivariant;
\item [(3)] $H$ is uniformly continuous in $t$, i.e. for any $\epsilon>0$, there exists $\delta>0$ such that 
\[d(H(x, t_1), H(x, t_2))\leq\epsilon\quad \text{for all}\ x\in\bigsqcup_{i\in I} X_i\text{ whenever}\ |t_1-t_2|<\delta;\]
\item [(4)] There exists a constant $C>0$, independent of $i\in I$ and $t\in [0, 1]$, such that
\[d(H(x,t),H(y,t))\leq C\cdot d(x,y)\quad\text{for all } x, y\in X_i,\ i\in I,\ t\in [0,1].\]
\end{enumerate}
Such a map $H$ is called an \emph{equivariant strong Lipschitz homotopy} between $g_0$ and $g_1$.
\end{definition}
\begin{definition}\label{def:collections-homotopy-equivalent}The collection $\{X_i\}_{i\in I}$ is said to be \emph{$F_i$-equivariant strongly Lipschitz homotopy equivalent} to $\{Y_i\}_{i\in I}$ if there exist $F_i$-equivariant uniformly Lipschitz maps 
\[g:\bigsqcup_{i\in I} X_i\to\bigsqcup_{i\in I}Y_i\quad\text{and}\quad h:\bigsqcup_{i\in I}Y_i\to\bigsqcup_{i\in I}X_i\]
satisfying:
\begin{enumerate}
\item [(1)] For each $i\in I$, $g(X_i)\subseteq Y_i$ and $h(Y_i)\subseteq X_i$;
\item [(2)] The composition $h\circ g$ is $F_i$-equivariant strongly Lipschitz homotopic to the identity map $\mathrm{id}_{\bigsqcup_{i\in I} X_i}$;
\item [(3)] The composition $g\circ h $ is $F_i$-equivariant strongly Lipschitz homotopic to the identity map $\mathrm{id}_{\bigsqcup_{i\in I} Y_i}$.
\end{enumerate}
\end{definition}
\begin{lemma}\label{lemma:simplex-evaluation-isomorphism}
Let $\{U_{i}\}_{i\in I}$ be a $\Gamma$-slice separated family as defined in Example \ref{example:slice-separated-sets}. For each $(d,R)\in\mathcal I$, let $\Delta_i^{(d)}(R)$ be the simplex in $P_{d}(X)$ spanned by the vertices 
$V_i(R):=\{v\in X\mid d_X(v,x_i)\leq R\}$.
Then the evaluation-at-zero map induces an isomorphism in $K$-theory:
\begin{equation*}e_{*}:K_{*}\bigl(\mathcal{A}^{*}_{L}(\Delta_i^{(d)}(R), i\in I)\bigr)
\xrightarrow{\cong}
K_{*}\bigl(\mathcal{A}^{*}(\Delta_i^{(d)}(R), i\in I)\bigr).
\end{equation*}
\end{lemma}
\begin{proof}[\textbf{Proof.}]\normalsize
For each $i\in I$, let
$\iota_i:\{x_i\}\hookrightarrow\Delta_i^{(d)}(R)$
be the natural inclusion and let
$r_i:\Delta_i^{(d)}(R)\to\{x_i\}$
be the constant map. Since $F_i$ fixes $x_i$, both maps are $F_i$-equivariant, and
$r_i\circ\iota_i=\operatorname{id}_{\{x_i\}}$.
Then
$$H_i(z,t)=(1-t)z+t x_i$$
is a well-defined homotopy on $\Delta_i^{(d)}(R)$, satisfies
$H_i(\cdot,0)=\operatorname{id}$,
$H_i(\cdot,1)=\iota_i\circ r_i$,
and is $F_i$-equivariant. Since $d, R$ are fixed and $X$ has bounded geometry, the simplices $\Delta_i^{(d)}(R)$ have uniformly bounded dimension and diameter. Hence,$H_i$ have a common Lipschitz constant and are uniformly
continuous in $t$, uniformly in $i\in I$. Therefore
$\{\Delta_i^{(d)}(R)\}_{i\in I}$ and $\{\{x_i\}\}_{i\in I}$ are
$F_i$-equivariantly strongly Lipschitz homotopy equivalent.

Since the families
$\iota=(\iota_i)_{i\in I}$,
$r=(r_i)_{i\in I}$,
and the homotopies $\{H_i\}_{i\in I}$ have control constants
independent of $i$, they induce well-defined homomorphisms on the
uniformly controlled product localization algebras.
By the uniform family version of \cite[Theorem 3.2]{Yu97}, the $K$-theory of localization algebras satisfies uniform $F_i$-equivariant strong Lipschitz homotopy invariance. The original proof applies if all control constants are taken uniformly in $i\in I$. Therefore, the induced maps $\iota_{L,*}$ and $r_{L,*}$ are mutually inverse isomorphisms. In particular,
\[\iota_{L,*} \colon
K_{*}\bigl(\mathcal{A}_L^*(\{x_i\},i\in I)\bigr)
\longrightarrow
K_{*}\bigl(\mathcal{A}_L^*(\Delta_i^{(d)}(R),i\in I)\bigr)\]
is an isomorphism.
The uniformly controlled $F_i$-equivariant coarse equivalences
$\{\iota_i\}_{i\in I}$ induce an isomorphism $\iota_{*}$ on the
$K$-theory of the corresponding Roe algebras. We obtain
the commutative diagram
\[
\begin{tikzcd}[row sep=large, column sep=huge]
K_{*}\bigl(\mathcal{A}_L^{*}(\{x_i\},i\in I)\bigr)
    \arrow[r,"e_{*}", "\cong"'] \arrow[d,"\iota_{L,*}"', "\cong"]
&
K_{*}\bigl(\mathcal{A}^{*}(\{x_i\},i\in I)\bigr)
    \arrow[d,"\iota_{*}", "\cong"']
\\
K_{*}\bigl(\mathcal{A}_L^{*}(\Delta_i^{(d)}(R),i\in I)\bigr)
    \arrow[r,"e_{*}"]
&
K_{*}\bigl(\mathcal{A}^{*}(\Delta_i^{(d)}(R),i\in I)\bigr).
\end{tikzcd}
\]
For the family $\{x_i\}_{i\in I}$, propagation is
zero. The standard Eilenberg--Swindle for localization
algebras shows that the kernel of evaluation-at-zero map has thetrivial
$K$-theory \cite{WY2020}. Hence the upper horizontal map is an isomorphism.
Since both vertical maps are isomorphisms, the lower horizontal
map is also an isomorphism. 
\end{proof}
\begin{proposition}\label{proposition:eval-slice-separated}
Suppose that the action $\Gamma\curvearrowright X$ has equivariant bounded
geometry and that $\mathcal{A}$ is a $(\Gamma,X)$-controlled proper
$C^{*}$-algebra. Let $\{U_i\}_{i\in I}$ be a $\Gamma$-slice
separated family as defined in
Example \ref{example:slice-separated-sets}, and set $U:=\bigcup_{i\in I}U_i$.
Then the evaluation-at-zero map induces an isomorphism in $K$-theory:
\[
e_{*}:\lim_{d\to\infty}K_{*}\bigl(C^{*}_{L}(P_{d}(X),\mathcal{A})^{\Gamma}_{U}\bigr)
\xrightarrow{\cong}
\lim_{d\to\infty}K_{*}\bigl(C^{*}(P_{d}(X),\mathcal{A})^{\Gamma}_{U}\bigr).
\]
\end{proposition}
\begin{proof}[\textbf{Proof.}]\normalsize
By Lemma \ref{lemma 3.16} (3) and (4), and by the continuity of $K$-theory with respect to inductive limits of $C^{*}$-algebras, we have the following commutative diagram:
\[
\begin{tikzcd}[row sep=large, column sep=large]
\displaystyle\lim_{d \to \infty } K_{*}\bigl(C^{*}_{L}(P_{d}(X),\mathcal{A})^{\Gamma}_{U}\bigr)
\arrow[r,"e_{*}"]
\arrow[d,"\cong"']
&
\displaystyle\lim_{d \to \infty } K_{*}\bigl(C^{*}(P_{d}(X),\mathcal{A})^{\Gamma}_{U}\bigr)
\arrow[d,"\cong"]
\\
\displaystyle\lim_{(d,R)\in\mathcal I}
K_{*}\bigl(\mathcal{A}^{*}_{L}(\Delta_i^{(d)}(R), i\in I)\bigr)
\arrow[r,"e_{*}"]
&
\displaystyle\lim_{(d,R)\in\mathcal I}
K_{*}\bigl(\mathcal{A}^{*}(\Delta_i^{(d)}(R), i\in I)\bigr).
\end{tikzcd}
\]
By Lemma \ref{lemma:simplex-evaluation-isomorphism}, the bottom horizontal map is an isomorphism for each fixed $(d,R)\in\mathcal I$, and hence remains an isomorphism after passing to the inductive limit over $\mathcal{I}$. Since the two vertical maps are isomorphisms, the top horizontal map is also an isomorphism. This completes the proof.
\end{proof}
In Proposition \ref{proposition:eval-slice-separated}, the finite subgroup $F_{i}=\Gamma_{x_i}$ fixes $x_{i}\in X$. To handle general $(\Gamma,X_0,r)$-separated families in Theorem \ref{thm:evaluation-map-isomorphism}, we need to consider the case where $F_i$ do not fix a vertex in $X$, but instead a point $q_i$ (such as the barycenter) in the Rips complex. We have the following result.
\begin{corollary}\label{prop:anchored-slice-evaluation}
Fix $D,s>0$. For every $i\in I$, let $q_i\in P_{d}(X)$ with support $E_i:=\operatorname{Supp}(q_i)\subseteq X$. Let
$F_i\leq\Gamma$ be a finite subgroup fixing $q_i$, and 
$V_{i,0}\subseteq M$ be an $F_i$-invariant open subset satisfying
\[V_{i,0}\subseteq O_{q_i,s}:=\bigcup_{v\in E_i}O_{v,s}.\]
Set
$V:=\bigsqcup_{i\in I}\bigsqcup_{\gamma F_i\in\Gamma/F_i}
\gamma V_{i,0}$.
Suppose that for every $S>0$,
$\sup_{x\in X}\#\left\{(i,\gamma F_i)\mid d_X(x,\gamma E_i)\leq S
\right\}<\infty$
and $\sup_i|F_i|<\infty$. Then evaluation-at-zero map
induces an isomorphism
\[e_{*}:\lim_{d\to\infty}
K_{*}\bigl(C_L^*(P_{d}(X),\mathcal{A})^\Gamma_V\bigr)
\longrightarrow
\lim_{d\to\infty}
K_{*}\bigl(C^*(P_{d}(X),\mathcal{A})^\Gamma_V\bigr).\]
\end{corollary}
\begin{proof}[\textbf{Proof.}]\normalsize
Since $F_i$ fixes $q_i$, it preserves $E_i=\operatorname{Supp}(q_i)$. Moreover,$\operatorname{diam}(E_i)\leq D$. For $d\geq D, R>0$, set
\[B_i^{(d)}(R):=\left\{z\in P_{d}(X)\mid d_X(v,E_i)\leq R
\ \text{for every}\ v\in\operatorname{Supp}(z)\right\}.\]
Replacing the balls centred at $x_i$ in Lemma \ref{lemma 3.16} by $B_i^{(d)}(R)$, the maps $T\mapsto(T_{e,i})_{i\in I}$ will give
\[C^{*}(P_{d}(X),\mathcal{A})^\Gamma_V\cong
\lim_{R\to\infty}\mathcal{A}^*\bigl(B_i^{(d)}(R),i\in I\bigr)\ \text{and}
\ C_L^*(P_{d}(X),\mathcal{A})^\Gamma_V\cong
\lim_{R\to\infty}
\mathcal{A}_L^*\bigl(B_i^{(d)}(R),i\in I\bigr).\]

We only check the control conditions. Suppose $T_{e,i}(z,w)\neq0$, and let $r_T$ be the support bound of $T$. Since
$V_{i,0}\subseteq O_{q_i,s}$, there are 
$v_z\in\operatorname{Supp}(z)$ and $v_i\in E_i$ such that
$O_{v_z,r_T}\cap O_{v_i,s}\neq\varnothing$.
Choose $S>0$ with
$\rho_-(S)\geq\max\{r_T,s\}$. By Definition \ref{def:proper-Gamma-X-algebra}(2)(b), we have that
$d_X(v_z,v_i)\leq2S$. Hence
$d_X(v,E_i)\leq2S+d$ for every
$v\in\operatorname{Supp}(z)$. Finite propagation gives the bound for $w$, uniformly in $i$.

Conversely, if $d_X(v,E_i)\leq R$, then
\[V_{i,0}\subseteq O_{v,\,s+\rho_+(R+D)}.\]
Thus the inverse map has a uniform support bound. The assumed bound on
\[\#\{(i,\gamma F_i)\mid d_X(x,\gamma E_i)\leq R\}\]
gives uniform local finiteness and uniform row and column bounds. 

For $d\geq2R+D$, the set
$\{v\in X\mid d_X(v,E_i)\leq R\}$ has diameter at most $2R+D$. So it spans a simplex
$\Delta_i^{(d)}(R)\subseteq P_{d}(X)$. In this case,
$B_i^{(d)}(R)=\Delta_i^{(d)}(R)$. Hence, for
$\mathcal I_{d}:=\{(d,R):d\geq2R+D\}$, we have that
\[\lim_{d\to\infty}C^*(P_{d}(X),\mathcal{A})^\Gamma_V
\cong\lim_{(d,R)\in\mathcal I_{d}}
\mathcal{A}^*\bigl(\Delta_i^{(d)}(R),i\in I\bigr).
\]The case of localization algebras is similar.

Each $\Delta_i^{(d)}(R)$ is $F_i$-invariant. The homotopy
$H_i(z,t)=(1-t)z+tq_i$ is $F_i$-equivariant, with control uniform in $i$. Lemma \ref{lemma:simplex-evaluation-isomorphism} therefore shows that evaluation-at-zero induces an isomorphism for every pair $(d,R)\in\mathcal I_{d}$. Passing to the direct limit proves the result.
\end{proof}
Now, we extend Proposition \ref{proposition:eval-slice-separated} to the case of 
$(\Gamma, X_0, r)$-separated family.
\begin{theorem}\label{thm:evaluation-map-isomorphism}
Let $\Gamma\curvearrowright X$ have equivariant bounded geometry, and let $\mathcal{A}$ be a $(\Gamma,X)$-controlled proper $C^{*}$-algebra with associated space $M$ as in
Definition \ref{def:proper-Gamma-X-algebra}. Let $r>0$, and let $X_0$ be a locally finite, $\Gamma$-invariant countable subset of $X$. Let $\{U_x\}_{x\in X_0}$ be a $(\Gamma,X_0,r)$-separated family as in Definition \ref{definition 3.10}, and set
$U=\bigcup_{x\in X_0}U_x$. Then the evaluation-at-zero map induces an isomorphism
\[e_{*}:\lim_{d\to\infty}K_{*}\bigl(
C_L^*(P_{d}(X),\mathcal{A})^\Gamma_U\bigr)
\xrightarrow{\cong}\lim_{d\to\infty}K_{*}\bigl(
C^*(P_{d}(X),\mathcal{A})^\Gamma_U\bigr).\]
\end{theorem}
\begin{proof}[\textbf{Proof.}]\normalsize 
Choose $R>0$ such that $\rho_-(R)\geq r$.
If $m\in U_x\cap U_y$, then
$O_{x,r}\cap O_{y,r}\neq\varnothing$.
Since $O_{x,r}\subseteq O_{x,\rho_-(R)}$ and similarly for $y$, Definition \ref{def:proper-Gamma-X-algebra}(2)(b) implies
$d_X(x,y)\leq 2R$. Denote 
\[\nu:=\sup_{x\in X}\#B_X(x,2R).
\]Then no point in $U$ belongs to more than $\nu$ sets $U_x$. Let $K\subseteq U$
be compact. Fix $y\in X$, Definition \ref{def:proper-Gamma-X-algebra}(2)(a) implies that there exists 
$t>0$ such that
$K\subseteq O_{y,t}$. Choose $S>0$ with
$\rho_-(S)\geq\max\{r,t\}$.
If $K\cap U_x\neq\varnothing$, then
$O_{y,t}\cap O_{x,r}\neq\varnothing$,
and hence $d_X(x,y)\leq2S$. Thus, as $X_0$ is locally finite, there are only finitely many such $x$.

Since $U$ is a second countable, locally compact proper
$\Gamma$-space, and $\{U_x\}_{x\in X_0}$ is a locally finite
$\Gamma$-invariant cover of $M$, there exists a $\Gamma$-equivariant
partition of unity
$\{\phi_x\}_{x\in X_0}$, satisfies 
$\operatorname{supp}(\phi_x)\subseteq U_x$ and 
$\phi_{\gamma x}(\gamma m)=\phi_x(m)$.
Let $\mathcal N$ be the nerve \cite{borsuk1948imbedding} of the cover. Then $\dim\mathcal N\leq\nu-1$.
The partition defines a continuous $\Gamma$-equivariant map
\[\Phi:U\longrightarrow|\mathcal N|,\quad
\Phi(m)=\sum_{x\in X_0}\phi_x(m)x.\]
Note that every simplex $\sigma$ in $\mathcal{N}$ spans a simplex in $P_{2R}(X)$, whose barycenter we denote by $b_\sigma$.

Let $\operatorname{star}(b_\sigma)$ be the open star of the
vertex $b_\sigma$ in the barycentric subdivision
$\operatorname{sd}\mathcal N$, and put
\[W_{\sigma}:=\Phi^{-1}\bigl(\operatorname{star}(b_\sigma)\bigr).
\]
Then $W_{\sigma}$ is open and
\[W_\sigma\subseteq\bigcap_{x\in\sigma}U_x
\subseteq O_{b_\sigma,r}.
\]
Moreover,
$\gamma W_\sigma=W_{\gamma\sigma}$.
For $0\leq q\leq\nu-1$, define
\[Y_q:=\bigsqcup_{\dim\sigma=q}W_\sigma.\]
It is a disjoint union since distinct vertices of the same dimension in $\operatorname{sd}\mathcal N$ have disjoint open stars. Each $Y_q$ is open and
$\Gamma$-invariant, and the open stars cover $|\mathcal N|$.
Therefore
\[U=\bigcup_{q=0}^{\nu-1}Y_q.\]\par
\textbf{Claim.}
For every nonempty $J\subseteq\{0,\ldots,\nu-1\}$, the set
$Y_J:=\bigcap_{q\in J}Y_q$ satisfies the assumptions of
Proposition \ref{prop:anchored-slice-evaluation}.
Set $m:=\max J$. Then $Y_J$ is the disjoint union of
$W_{\boldsymbol\sigma}:=\bigcap_{q\in J}W_{\sigma_q}$, where
$\boldsymbol\sigma=(\sigma_q)_{q\in J}$ runs through chains with $\dim\sigma_q=q$. Indeed, if two such intersections meet, all their faces are comparable \cite{Mun84}. Faces of the same dimension must then be equal.
Define
\[F_{\boldsymbol\sigma}:=\{\gamma\in\Gamma\mid
\gamma\sigma_q=\sigma_q\text{for all}q\in J\},
\quad q_{\boldsymbol\sigma}=b_{\sigma_m}\in P_{2R}(X).
\]
Then $F_{\boldsymbol\sigma}$ fixes $q_{\boldsymbol\sigma}$.
Its action on the vertices of $\sigma_m$ has degree at most $\nu$,
and its kernel is contained in the stabilizer of every vertex. Hence,
\[|F_{\boldsymbol\sigma}|\leq\nu!\sup_{x\in X}|\Gamma_x|<\infty.
\]
Moreover,
$W_{\boldsymbol\sigma}\subseteq W_{\sigma_m}
\subseteq O_{q_{\boldsymbol\sigma},r}$, and
if $\gamma W_{\boldsymbol\sigma}\cap W_{\boldsymbol\sigma}\neq\varnothing$, then $\gamma\boldsymbol\sigma=\boldsymbol\sigma$. This means that $\gamma\in F_{\boldsymbol\sigma}$.
Thus, if we choose one chain from each $\Gamma$-orbit,
$Y_J$ is a disjoint union of slices
$\Gamma\times_{F_{\boldsymbol\sigma}}W_{\boldsymbol\sigma}$.

Let
$E_{\boldsymbol\sigma}$ be the vertex set of $\sigma_m$. Then
$|E_{\boldsymbol\sigma}|\leq\nu$ and
$\operatorname{diam}(E_{\boldsymbol\sigma})\leq2R$. If
$d_X(x,\gamma E_{\boldsymbol\sigma})\leq S$, 
$$\gamma E_{\boldsymbol\sigma}\subseteq B_X(x,S+2R).$$ For fixed $J$, each maximal face belongs to only
uniformly many chains. We have that
\[\sup_{x\in X}\#\left\{(\boldsymbol\sigma,\gamma F_{\boldsymbol\sigma}):
d_X(x,\gamma E_{\boldsymbol\sigma})\leq S
\right\}<\infty.\]

Let $D=2R$ and
$s=r$ in
Corollary \ref{prop:anchored-slice-evaluation}, we obtian
\[e_{*}:\lim_{d\to\infty}K_{*}\bigl(C_L^*(P_{d}(X),\mathcal{A})^\Gamma_{Y_J}\bigr)\xrightarrow{\cong}\lim_{d\to\infty}
K_{*}\bigl(C^*(P_{d}(X),\mathcal{A})^\Gamma_{Y_J}\bigr).\]
Using a Mayer-Vietoris argument to glue the pieces together from the intersections $Y_J$, the result for $U$ follows.
\end{proof}

\section{Proof of the Main Theorem}

We now turn to the proof of the main theorem, beginning with the following lemma.
\begin{lemma}\label{lemma 2.25}
Let $(X,d_X)$ be a metric space of bounded geometry, and let
$\Gamma\curvearrowright X$ be an isometric action. Let
$X_0\subseteq X$ be a $\Gamma$-invariant subset. For $R>0$, define a
graph $\mathcal{G}_{R}$ whose vertices are the $\Gamma$-orbits in
$X_0$. Two distinct vertices $\mathcal{O}_1,\mathcal{O}_2\in X_{0}/\Gamma$
are joined by an edge if there exist $x_1\in\mathcal{O}_1$ and
$x_2\in\mathcal{O}_2$ such that
$d(x_1,x_2)\leq 2R$.
Then $\mathcal{G}_R$ has uniformly bounded degree.
\end{lemma}
\begin{proof}[\textbf{Proof.}]\normalsize 
Let $\mathcal{O}_1\in X_0/\Gamma$ be a given vertex, and fix a representative $x_1\in\mathcal{O}_1$. Suppose $\mathcal{O}_2$ is a neighbor of $\mathcal{O}_1$ in $\mathcal{G}_R$. By the definition of the edges in $\mathcal{G}_R$, there exist $y_1\in\mathcal{O}_1$ and $y_2\in\mathcal{O}_2$ such that $d(y_1, y_2)\leq 2R$. Since $y_1\in\mathcal{O}_1$, there exists $\gamma\in \Gamma$ such that $y_1=\gamma x_1$. Then
\[d(x_1,\gamma^{-1}y_2)=d(\gamma x_1, y_2)=d(y_1, y_2)\leq 2R.\]
Since $\gamma^{-1}y_2$ belongs to the same orbit $\mathcal{O}_2$, every neighboring orbit $\mathcal{O}_2$ must have a representative in the ball $B(x_1,2R)$. Since $X_0$ is locally finite, the number of points in $X_0\cap B(x_1,2R)$ is uniformly bounded by a constant $C_{2R}$, independent of $x_1$. Thus, the number of distinct $\Gamma$-orbits intersecting $B(x_1,2R)$ is at most $C_{2R}$. Hence, the degree of any vertex in $\mathcal{G}_R$ is uniformly bounded by $C_{2R}$.
\end{proof}
\begin{proof}[\textbf{Proof of Theorem \ref{main theorem}.}]\normalsize
Recall from Definition \ref{def:proper-Gamma-X-algebra} that for each $x\in X$ and $r>0$, we associate an open subset $O_{x,r}\subseteq M$ satisfying conditions (a)--(d). Since $\Gamma\curvearrowright X$ has equivariant bounded geometry, $X$ is a countable, locally finite, and $\Gamma$-invariant index set. For simplicity, we take $X_{0}:=X$ as our index set and define
\[O_{r}:=\bigcup_{x\in X_0} O_{x,r}\subseteq M.\]
Fix $d>0$. For $r'>r>0$ we have inclusions
$C^{*}(P_{d}(X),\mathcal{A})^\Gamma_{O_r}\subseteq
C^{*}(P_{d}(X),\mathcal{A})^\Gamma_{O_{r'}}$. We claim these ideals exhaust the full algebras. Indeed, let
$T\in\mathbb{C}[P_{d}(X),\mathcal{A}]^\Gamma$. By
Definition \ref{def:equiv-alg-roe-X}(5) there is a single $r_T>0$ with
$\operatorname{Supp}(T(x,y))\subseteq O_{x,r_T}$ for all $x,y$. Since
$X_0=X$, the coefficient support of $T$ satisfies
\[\bigcup_{x,y}\operatorname{Supp}(T(x,y))\subseteq
\bigcup_{x\in X_0}O_{x,r_T}=O_{r_T},\]
so $T\in\mathbb{C}[P_{d}(X),\mathcal{A}]^\Gamma_{O_{r_T}}$. The same argument,
using the $t$-independent bound in Definition \ref{def:equiv-alg-loc-X}(2),
applies to $\mathbb{C}_{L}[P_{d}(X),\mathcal{A}]^\Gamma$. Hence,
$\bigcup_{r>0}\mathbb{C}[P_{d}(X),\mathcal{A}]^\Gamma_{O_r}
=\mathbb{C}[P_{d}(X),\mathcal{A}]^\Gamma$
(and likewise for the localization algebra), and passing to norm completions
we have that
\[\lim_{r\to\infty }C^{*}(P_{d}(X),\mathcal{A})^\Gamma_{O_r}
=C^{*}(P_{d}(X),\mathcal{A})^\Gamma
\quad\text{and}\quad
\lim_{r\to\infty}C_L^{*}(P_{d}(X),\mathcal{A})^\Gamma_{O_r}
=C_{L}^{*}(P_{d}(X),\mathcal{A})^\Gamma.\]

For a fixed $r>0$, if $d'>d>0$, the inclusions $P_{d}(X)\hookrightarrow P_{d'}(X)$ induce:
\[C^{*}_{L}(P_{d}(X),\mathcal{A})^{\Gamma}_{O_{r}}\hookrightarrow C^{*}_{L}(P_{d'}(X),\mathcal{A})^{\Gamma}_{O_{r}}\quad \text{and}\quad C^{*}(P_{d}(X),\mathcal{A})^{\Gamma}_{O_{r}}\hookrightarrow C^{*}(P_{d'}(X),\mathcal{A})^{\Gamma}_{O_{r}}.\]
We have the following commutative diagram:
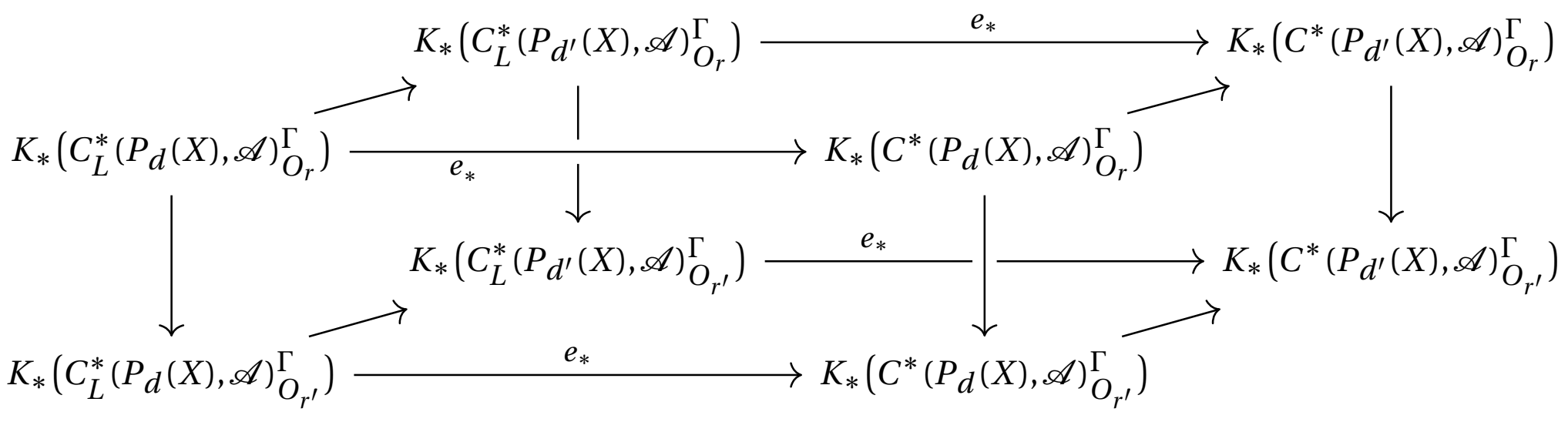
\begin{center}
\small
\begin{tikzcd}[row sep=small, column sep=small]
& K_{*}\bigl(C^{*}_{L}(P_{d'}(X),\mathcal{A})_{O_{r}}^{\Gamma}\bigr) \arrow[dd] \arrow[rr,"e_{*}"] & & K_{*}\bigl(C^{*}(P_{d'}(X),\mathcal{A})_{O_{r}}^{\Gamma}\bigr) \arrow[dd] \\
K_{*}\bigl(C^{*}_{L}(P_{d}(X),\mathcal{A})_{O_{r}}^{\Gamma}\bigr) \arrow[ru] \arrow[dd] \arrow[rr, crossing over, "e_{*}^{\ }"' near start] & & K_{*}\bigl(C^{*}(P_{d}(X),\mathcal{A})_{O_{r}}^{\Gamma}\bigr) \arrow[ru] & \\
& K_{*}\bigl(C^{*}_{L}(P_{d'}(X),\mathcal{A})_{O_{r'}}^{\Gamma}\bigr) \arrow[rr,"e_{*}" near start] & & K_{*}\bigl(C^{*}(P_{d'}(X),\mathcal{A})_{O_{r'}}^{\Gamma}\bigr) \\
K_{*}\bigl(C^{*}_{L}(P_{d}(X),\mathcal{A})_{O_{r'}}^{\Gamma}\bigr) \arrow[ru] \arrow[rr,"e_{*}"] & & K_{*}\bigl(C^{*}(P_{d}(X),\mathcal{A})_{O_{r'}}^{\Gamma}\bigr) \arrow[ru] \arrow[from=uu, crossing over] &
\end{tikzcd}
\end{center}
Since the $K$-theory functor commutes with direct limits, we can interchange the order of limits as shown in the following commutative square:
\begin{center}
\begin{tikzcd}[row sep=large, column sep=large]
{\lim\limits_{d\rightarrow\infty}K_{*}(C^{*}_{L}(P_{d}(X),\mathcal{A})^{\Gamma}) }\arrow[r,"e_{*}^{\A}"] \arrow[d,"{\cong}"] & {\lim\limits_{d\rightarrow\infty}K_{*}(C^{*}(P_{d}(X),\mathcal{A})^{\Gamma})} \arrow[d,"{\cong}"] \\
{\lim\limits_{d\rightarrow\infty}\lim\limits_{r\rightarrow\infty}K_{*}(C^{*}_{L}(P_{d}(X),\mathcal{A})_{O_{r}}^{\Gamma})} \arrow[d,"\cong"] \arrow[r,"e_{*}"] & {\lim\limits_{d\rightarrow\infty}\lim\limits_{r\rightarrow\infty}K_{*}(C^{*}(P_{d}(X),\mathcal{A})_{O_{r}}^{\Gamma})}\arrow[d,"\cong"] \\ {\lim\limits_{r\rightarrow\infty}\lim\limits_{d\rightarrow\infty}K_{*}(C^{*}_{L}(P_{d}(X),\mathcal{A})_{O_{r}}^{\Gamma})}  \arrow[r,"e_{*}"] & {\lim\limits_{r\rightarrow\infty}\lim\limits_{d\rightarrow\infty}K_{*}(C^{*}(P_{d}(X),\mathcal{A})_{O_{r}}^{\Gamma})}
\end{tikzcd}
\end{center}
Thus, to prove that the top horizontal map $e_{*}^{\A}$ is an isomorphism, it suffices to prove that for each fixed $r>0$, the evaluation map
\begin{equation}\label{eq:fixed-r-isomorphism}
e_{*}:\lim_{d\to\infty}K_{*}\bigl(C^{*}_{L}(P_{d}(X),\mathcal{A})_{O_{r}}^{\Gamma}\bigr)\rightarrow\lim_{d\to\infty}K_{*}\bigl(C^{*}(P_{d}(X),\mathcal{A})_{O_{r}}^{\Gamma}\bigr)
\end{equation}
is an isomorphism.

We now correspond the index set $X_{0}$ with a finite union of separated families. First, choose an integer $R$ large enough such that $\rho _{-}(R)>r$. We construct the decomposition inductively. Choose $x_1 \in X_{0}$ and let $X_1^{(1)}=\Gamma\cdot x_1$. Next, select $x_2\in X_0\setminus\bigcup_{\gamma\in\Gamma} B(\gamma\cdot x_1, 2R)$ and let $X_{1}^{(2)}=\Gamma\cdot x_{2}$. At step $k$, choose
\[x_k\in X_0\setminus\bigcup_{i=1}^{k-1}\bigcup_{\gamma \in\Gamma}B(\gamma\cdot x_i,2R),\]
and define $X_1^{(k)}=\Gamma\cdot x_k$. Let $X_1=\bigsqcup_{k}X_1^{(k)}$. We repeat this greedy algorithm \cite{BM08} on $X_{0}\setminus X_1$ to extract $X_2$, and so on. 
This procedure is equivalent to extracting independent sets from the orbit graph $\mathcal{G}_R$ on $X_0/\Gamma$, where edges connect orbits at a distance $\leq 2R$. By Lemma \ref{lemma 2.25}, this graph $\mathcal{G}_R$ has a uniformly bounded degree, say $\Delta$. According to standard graph theory (Cf. \cite{Die17}), any graph with maximum degree $\Delta$ can be colored with at most $N:=\Delta+1$ colors. Consequently, the greedy algorithm is guaranteed to terminate and partitions $X_{0}$ into at most $N$ independent sets, yielding a finite disjoint orbit decomposition:
\[X_{0}=\bigsqcup_{n=1}^{N}X_n.\]

For each color class $n$, the distance between any two distinct orbits in $X_n$ is greater than $2R$. Since $r<\rho_{-}(R)$, Definition \ref{def:proper-Gamma-X-algebra}(2)(b) guarantees that $O_{\gamma x, r}\cap O_{\gamma' y,r}=\emptyset$ for distinct orbits. Thus, for each $n$, the family $\{O_{x,r}\}_{x\in X_n}$ forms a $(\Gamma, X_n,r)$-separated family of open sets in $M$.

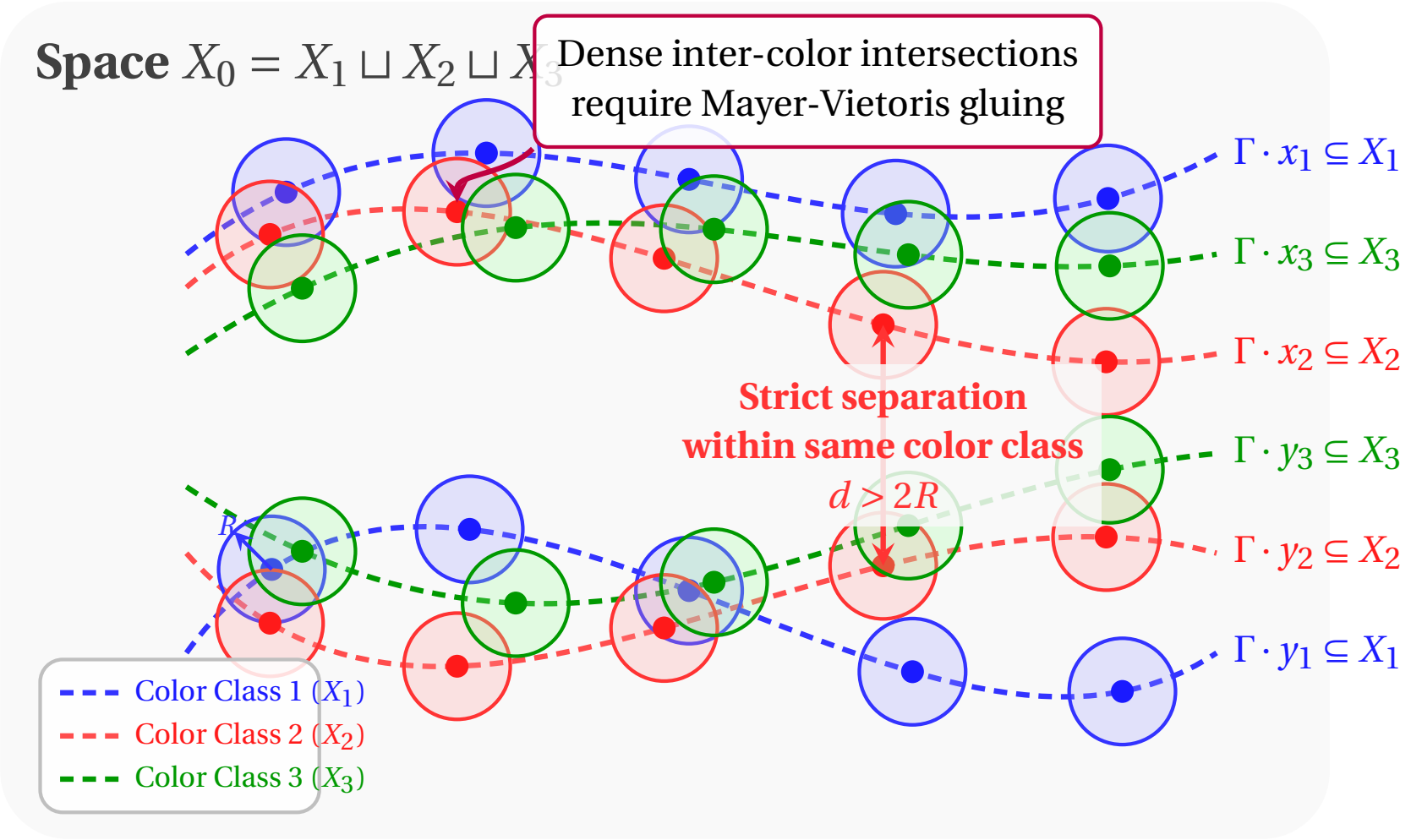
\begin{figure}[htbp]
\centering
\begin{tikzpicture}[
    scale=0.52,
    >=stealth,
    group1orb/.style={very thick, densely dashed, blue!80},
    group2orb/.style={very thick, densely dashed, red!70},
    group3orb/.style={very thick, densely dashed, green!60!black},
    group1ball/.style={fill=blue!40, fill opacity=0.25, draw=blue!80, thick},
    group2ball/.style={fill=red!40, fill opacity=0.25, draw=red!80, thick},
    group3ball/.style={fill=green!40, fill opacity=0.25, draw=green!60!black, thick},
    point/.style={circle, fill, inner sep=1.8pt}
]

\fill[gray!5, rounded corners=15pt] (-3.8, -0.8) rectangle (16.2, 11.8);

\node[font=\large\bfseries, darkgray, anchor=north west] at (-3.5, 11.4) {Space $X_{0} = X_1 \sqcup X_2 \sqcup X_3$};

\def\myrad{0.8}

\draw[group1orb] (-1, 8.0) 
    .. controls (4, 12.0) and (9, 6.5) .. (14.5, 9.5)
    node[pos=0.1, coordinate] (B1_1) {}
    node[pos=0.3, coordinate] (B1_2) {}
    node[pos=0.5, coordinate] (B1_3) {}
    node[pos=0.7, coordinate] (B1_4) {}
    node[pos=0.9, coordinate] (B1_5) {}
    node[right, font=\small, text=blue!90] {$\Gamma \cdot x_1 \subseteq X_1$};

\draw[group1orb] (-1, 2.0) 
    .. controls (3, 7.5) and (10, -1.0) .. (14.5, 2.0)
    node[pos=0.1, coordinate] (B2_1) {}
    node[pos=0.3, coordinate] (B2_2) {}
    node[pos=0.5, coordinate] (B2_3) {}
    node[pos=0.7, coordinate] (B2_4) {}
    node[pos=0.9, coordinate] (B2_5) {}
    node[right, font=\small, text=blue!90] {$\Gamma \cdot y_1 \subseteq X_1$};

\draw[group2orb] (-1, 7.5) 
    .. controls (3, 11.0) and (9, 5.5) .. (14.5, 6.5)
    node[pos=0.1, coordinate] (R1_1) {}
    node[pos=0.3, coordinate] (R1_2) {}
    node[pos=0.5, coordinate] (R1_3) {}
    node[pos=0.7, coordinate] (R1_4) {}
    node[pos=0.9, coordinate] (R1_5) {}
    node[right, font=\small, text=red!90] {$\Gamma \cdot x_2 \subseteq X_2$};

\draw[group2orb] (-1, 3.5) 
    .. controls (3, -1.0) and (9, 5.0) .. (14.5, 3.5)
    node[pos=0.1, coordinate] (R2_1) {}
    node[pos=0.3, coordinate] (R2_2) {}
    node[pos=0.5, coordinate] (R2_3) {}
    node[pos=0.7, coordinate] (R2_4) {}
    node[pos=0.9, coordinate] (R2_5) {}
    node[right, font=\small, text=red!90] {$\Gamma \cdot y_2 \subseteq X_2$};

\draw[group3orb] (-1, 6.5) 
    .. controls (5, 10.5) and (9, 7.0) .. (14.5, 8.0)
    node[pos=0.1, coordinate] (G1_1) {}
    node[pos=0.3, coordinate] (G1_2) {}
    node[pos=0.5, coordinate] (G1_3) {}
    node[pos=0.7, coordinate] (G1_4) {}
    node[pos=0.9, coordinate] (G1_5) {}
    node[right, font=\small, text=green!60!black] {$\Gamma \cdot x_3 \subseteq X_3$};

\draw[group3orb] (-1, 4.5) 
    .. controls (5, 0.5) and (9, 4.5) .. (14.5, 5.0)
    node[pos=0.1, coordinate] (G2_1) {}
    node[pos=0.3, coordinate] (G2_2) {}
    node[pos=0.5, coordinate] (G2_3) {}
    node[pos=0.7, coordinate] (G2_4) {}
    node[pos=0.9, coordinate] (G2_5) {}
    node[right, font=\small, text=green!60!black] {$\Gamma \cdot y_3 \subseteq X_3$};

\foreach \i in {1,2,3,4,5} {
    \draw[group1ball] (B1_\i) circle (\myrad); \node[point, blue!90] at (B1_\i) {};
    \draw[group1ball] (B2_\i) circle (\myrad); \node[point, blue!90] at (B2_\i) {};
    \draw[group2ball] (R1_\i) circle (\myrad); \node[point, red!90] at (R1_\i) {};
    \draw[group2ball] (R2_\i) circle (\myrad); \node[point, red!90] at (R2_\i) {};
    \draw[group3ball] (G1_\i) circle (\myrad); \node[point, green!60!black] at (G1_\i) {};
    \draw[group3ball] (G2_\i) circle (\myrad); \node[point, green!60!black] at (G2_\i) {};
}

\draw[<->, very thick, red!80] (R1_4) -- node[fill=gray!5, fill opacity=0.85, text opacity=1, inner sep=4pt, font=\small\bfseries, align=center] {Strict separation\\within same color class\\$d > 2R$} (R2_4);

\draw[->, thick, blue!80] (B2_1) -- node[above left, font=\footnotesize] {$R$} ++(135:\myrad);

\node[fill=white, fill opacity=0.85, text opacity=1, draw=purple, thick, rounded corners=3pt, inner sep=5pt, font=\small, align=center] (overlapbox) at (8.5, 10.6) {Dense inter-color intersections\\require Mayer-Vietoris gluing};
\draw[->, very thick, purple, shorten >=2pt] (overlapbox.south west) to[out=-135, in=90] (R1_2);

\begin{scope}[shift={(-3.2, -0.4)}]
    \fill[white, fill opacity=0.85, text opacity=1, draw=gray!50, rounded corners=5pt, thick] (0, 0) rectangle (4.2, 2.3);
    
    \draw[group1orb] (0.3, 1.8) -- (1.2, 1.8);
    \node[right, font=\footnotesize, text=blue!90] at (1.2, 1.8) {Color Class 1 ($X_1$)};
    
    \draw[group2orb] (0.3, 1.15) -- (1.2, 1.15);
    \node[right, font=\footnotesize, text=red!90] at (1.2, 1.15) {Color Class 2 ($X_2$)};
    
    \draw[group3orb] (0.3, 0.5) -- (1.2, 0.5);
    \node[right, font=\footnotesize, text=green!60!black] at (1.2, 0.5) {Color Class 3 ($X_3$)};
\end{scope}

\end{tikzpicture}
\caption{Visualization of the orbit decomposition via graph coloring. To overcome the non-cocompactness of the space, the entire index set $X_0$ is divided into a finite number of color groups (represented here by 3 colors). Orbits belonging to different color classes can closely overlap and cross each other, creating heavily intersecting $C^{*}$-ideals that are linked together via Mayer-Vietoris exact sequences. Conversely, orbits within the same color group (e.g., the two distinct red orbits) maintain a strict spacing rule of $d > 2R$. This ensures that each individual color class forms a completely separate, non-overlapping set, allowing local structures to extend to the overall space.}
\label{fig:orbit-decomposition-3colors}
\end{figure}

For each $n=1,\dots, N$, we define:
\[O_{n,r}:=\bigcup_{x\in X_n}O_{x,r}.\]
Since $\{O_{x,r}\}_{x\in X_n}$ is a $(\Gamma,X_n,r)$-separated family, Theorem \ref{thm:evaluation-map-isomorphism} directly guarantees that the evaluation map is an isomorphism:
\[e_{*}^{(n)}:\lim_{d\to\infty}K_{*}\bigl(C_L^{*}(P_{d}(X),\mathcal{A})^\Gamma_{O_{n,r}}\bigr) 
\xrightarrow{\cong}\lim_{d\to\infty }K_{*}\bigl(C^{*}(P_{d}(X),\mathcal{A})^\Gamma_{O_{n,r}}\bigr).
\]

Finally, we glue these $N$ pieces together. Note that $O_{r}=\bigcup_{n=1}^{N}O_{n,r}$ is a finite open cover. For any subset of indices $I\subseteq\{1,\dots,N\}$, the intersection $\bigcap_{i\in I}O_{i,r}$ is an open $\Gamma$-invariant subset of $O_{k,r}$. Theorem \ref{thm:evaluation-map-isomorphism} guarantees that $e_{*}$ is an isomorphism for all such finite intersections. 

We proceed by induction on $N$ using the Mayer-Vietoris sequence. The case $N=1$ is exactly Theorem \ref{thm:evaluation-map-isomorphism}. For $N>1$, we write $O_r= O_{1,r}\cup\bigl(\bigcup_{n=2}^{N}O_{n,r}\bigr)$. By the induction hypothesis, $e_{*}$ is an isomorphism for $O_{1,r}$, for $\bigcup_{n=2}^{N}O_{n,r}$, and for their intersection. Applying the Five Lemma to the corresponding Mayer-Vietoris sequences in $K$-theory implies that $e_{*}$ is an isomorphism for the entire union $O_r$. 
This establishes the isomorphism for \ref{eq:fixed-r-isomorphism}, thereby completing the proof of Theorem \ref{main theorem}.
\end{proof}

\noindent\textbf{Funding.}
This work was supported by the NSFC (Nos. 12571135, 12171156), the
Science and Technology Commission of Shanghai (No. 22DZ2229014), the Key
Laboratory of MEA, Ministry of Education, and the Shanghai Key Laboratory of PMMP,
East China Normal University.
\vspace{\baselineskip}

\noindent\textbf{Declaration of AI Use}
The authors used ChatGPT to assist with English language editing and improving the clarity and presentation of the manuscript. The authors take full responsibility for the final content.

\bibliographystyle{alpha}
\bibliography{reff}
\printendmatterinfo
\end{document}